# Finite Dominating Sets for the Refueling Station Location Problem in Fleet Operations

Moddassir Khan Nayeem [a], Fuhad Ahmed Opu [a], Omar Abbaas [a*], and Sara Abu-Aridah [a]

[a] Department of Mechanical, Aerospace, and Industrial Engineering, The University of Texas at San Antonio, San Antonio, TX 78249, USA

## Abstract

This study considers a set of routes used by public transportation vehicles and dedicated distribution fleets in a general network. We aim to optimally locate alternative fuel refueling stations in the network to serve these dedicated routes. Deviations from prescribed routes for refueling purposes are allowed. Unlike most related literature, our approach considers all points in the network as candidate refueling station locations. We derive coverage constraints for any candidate location to serve a given route. Then we develop an exact algorithm to establish a finite dominating set (FDS) of candidate locations guaranteed to include an optimal solution to the problem. This set can be used in a mathematical model to minimize the number of stations required to cover all flows in the network. Numerical experiments on realistic networks are presented to illustrate the proposed methodology and to demonstrate its scalability and sensitivity to changes in parameter values.

(*) Corresponding author: Tel.: +1(210)458-8746; E-mail address: omar.abbaas@utsa.edu.

## 1. Introduction

Over the past 50 years, fossil fuels have supplied over 90% of the energy used in transportation, leading to significant Greenhouse Gas (GHG) emissions and contributing to climate change [1, 2]. In 2022, transportation accounted for about 28% of total GHG emissions in the USA, surpassing all other economy sectors [1]. Growing environmental concerns and energy independence efforts have accelerated the shift toward cleaner energy sources like electricity, hydrogen, and biofuels. Municipal fleet data shows that Alternative Fuel Vehicles (AFVs) can cut GHG emissions by up to 36% while maintaining operational efficiency, making them viable for long-term environmental sustainability goals [2]. Public transportation authorities in many major cities are adding electric vehicles to their fleets. For example, the Metropolitan Transportation Authority in New York City pledged to convert all 5,800 of its bus fleet to zero-emission vehicles by 2040 with plans to stop purchasing non-electric buses in 2029 [3]. Similarly, many logistics and parcel delivery fleets are transitioning to AFVs due to their economic and environmental benefits. For instance, in 2019 Amazon and Rivian announced plans to deploy 100,000 electric delivery vans by 2030 to more than 100 cities in the US [4].

The integration of Electric Vehicles (EVs) into commercial fleets poses logistical challenges due to their relatively limited driving range, long charging times, and inadequate infrastructure. Today, most EV models offer driving ranges between 100 and 300 miles per charge, but this can vary significantly based on factors such as vehicle type, battery capacity, and driving conditions [5]. Due to their limited driving range, EVs may need to recharge multiple times during the day which can make frequent trips to charging stations inefficient and impractical. To address the logistical challenges, a variety of strategies are being considered. One promising approach is the use of battery swapping stations, which allows depleted batteries to be exchanged quickly with fully charged ones [6]. Battery swapping has the potential to significantly reduce "refueling" times and extend the effective range of EVs [7]. However, widespread adoption of battery swapping stations faces several hurdles, including the need for standardization across vehicle manufacturers and significant infrastructure investment.

As infrastructure development is both costly and time-consuming, it is important to establish refueling stations at locations that cover the maximum traffic flow within a transportation network [8-12]. In traditional network location theory, the primary assumption is that demands originate at network nodes. This assumption makes the problem easier since network nodes form a natural finite set of candidate locations and are usually more accessible compared to points along the edges. However, people typically do not make trips to refueling stations solely for the purpose of refueling [13]. As a result, recent literature on refueling station locations usually assumes that demand occurs along preplanned trips. Upchurch and Kuby [14] present a comparison between the two approaches, namely node-based and flow-based models.

Flow-based models represent traffic flow using an Origin-Destination (OD) matrix, where each element indicates the number of vehicles traveling from the origin to the destination. The objective is to locate facilities to maximize the covered flows. Hodgson et al. [15] and Berman et al. [16] pioneered flow-based models and developed the flow-capturing location model (FCLM). Their FCLM determines the optimal placement of a given number of facilities to maximize the number of customers who encounter at least one facility along their predetermined path. By considering the flow of trips rather than focusing on node-based demands, flow-based models provide a more realistic representation of the refueling station location problem and enable the development of strategies that better serve the needs of AFVs.

While these models provide a more accurate representation of the refueling station location problem, a critical consideration is the transition from the discrete to the continuous version of the problem. The discrete version, commonly addressed in literature, limits candidate locations to a predetermined finite set, potentially leading to sub-optimal solutions. Conversely, in the continuous version, refueling stations can be placed anywhere in the network making the initial set of candidate locations infinite. This approach allows a comprehensive exploration of the solution space, increasing the likelihood of finding an optimal solution. By contrast, the discrete version limits the search to a predefined set of candidate locations, which makes the problem easier to solve but may exclude the optimal solution if it does not lie in this predetermined set.

This study uses the continuous approach to locate a set of refueling stations in a general network to serve vehicles dedicated to predetermined closed routes, such as public transportation fleets. We consider common AFVs' challenges such as limited driving range and underdeveloped refueling infrastructure leading to deviations from predetermined routes for refueling purposes. These deviations allow vehicles from different non-overlapping routes to share the use of refueling stations. Longer deviation distances enable vehicles to access refueling stations far from their routes which reduces the required number of refueling stations for covering all flows in the network. However, we limit the maximum allowed deviation distance for several practical reasons. First, limiting deviation distances limits the portion of the fuel tank/battery used for traveling to refueling stations, and consequently reduces cost. Second, it minimizes the inconvenience and wasted time caused by lengthy refueling trips. Finally, it limits the environmental impact of refueling trips. Based on these considerations, we develop the theoretical basis and algorithm to extract a Finite Dominating Set (FDS) of candidate locations in the traffic network. This set is guaranteed to include an optimal solution to the problem. The objective is to minimize the number of refueling stations required to cover all flows in the network. Our methodology allows us to find all optimal solutions, which gives decision-makers the flexibility to consider additional practical factors such as land prices, accessibility, etc. in their final decision. To the best of our knowledge, this paper is the first to use the continuous approach for the alternative fuel refueling station location problem on a general network for vehicles dedicated to closed routes.

The rest of this paper is organized as follows. Section 2 reviews the related literature. In Section 3, the problem statement is provided along with the assumptions and notation used throughout the paper. Section 4 presents the proposed methodology. In Section 5, we conduct a series of numerical experiments to illustrate the proposed methodology. We conclude with final remarks and suggestions for future research in Section 6.

## 2. Literature Review

The problem of optimally locating service facilities on traffic networks has attracted significant attention in the literature. Research in this area has evolved over the years addressing factors such as flow coverage, deviation from shortest path, multi-period, driving range, continuous and discrete approaches, as well as different types of networks. Figure 1 provides a visual summary of key models and their advancements. In Subsection 2.1, we review existing literature related to refueling station location models, summarizing key contributions and advancements in the field. Based on this review, Subsection 2.2 discusses the research gaps and highlights the uniquness of this study.

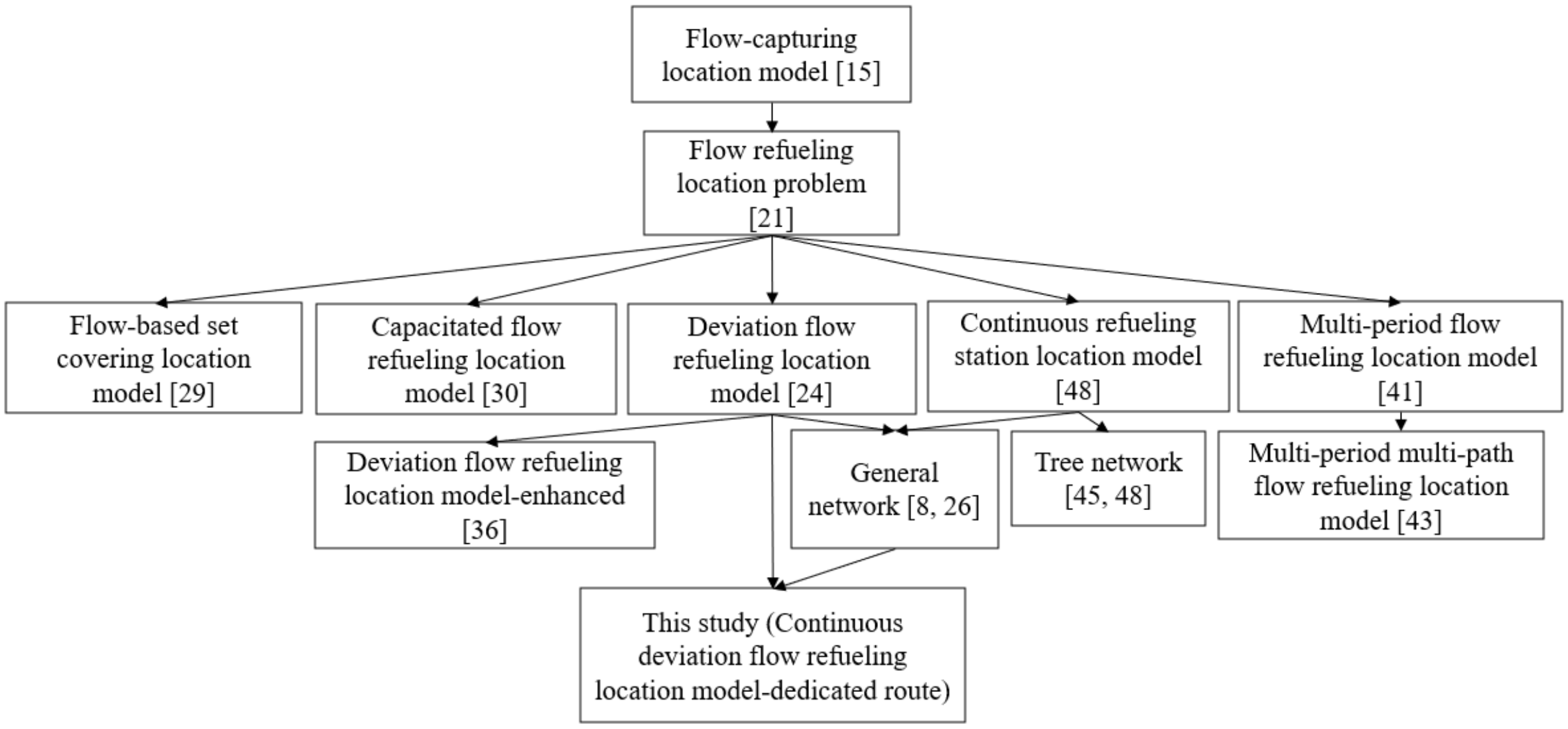

Figure 1: Evolution of AFV refueling station location models

### *2.1. Refueling Station Location Models*

Traditionally, facility location models have presumed that demand is generated at or attracted by service facilities [17-19]. While this assumption holds for certain facilities, such as restaurants and major retail stores, refueling stations operate differently; they typically serve as stops along preplanned routes rather than being the final destination. An early attempt to address this unique aspect is the FCLM proposed by Hodgson [15]. This model aims to locate facilities to maximize the number of customers who encounter at least one facility along their designated paths. Building upon this foundation, Berman et al. [16] introduced a model focused on identifying optimal locations for multiple service facilities, balancing two main objectives: maximizing the captured flow and minimizing the number of facilities needed to serve a predetermined portion of the total traffic flow.

The FCLM model considers a flow to be covered if it encounters at least one facility along its path. However, FCLM model does not consider the driving range which is a critical factor for AFVs. In reality, a traffic flow is only refueled if it can complete a round trip from origin to destination and back without running out of fuel. AFVs have a limited driving range per refill, requiring multiple refueling stops along their route [20]. To address this limitation, Kuby and Lim [21] introduced the flow refueling location problem (FRLP), which considers limited driving range and aims to establish a number of refueling stations on the network to maximize the covered round-trip travels. However, finding all combinations of refueling station locations that can cover all round trips in the network can make the model computationally intractable. To overcome this shortcoming, Lim and Kuby [22] proposed three heuristic algorithms: greedy-adding, greedy-adding with substitution, and genetic algorithm. Additionally, the FRLP model developed by Kuby and Lim [21] did not consider the coverage of path segments, which was later reformulated by Capar and Kuby [9] and Capar et al. [23].

The studies discussed above primarily focus on two key assumptions: first, customers receive services from alternative fuel stations located along their pre-planned shortest path, and second, alternative fuel stations have adequate capacity to serve all passing flows. However, during the initial transition period towards AFVs, the scarcity of alternative fuel stations may require drivers to deviate from their preplanned shortest paths to access refueling services [24-26]. Thus, FRLP can be more realistic by allowing deviations from preplanned paths for refueling purposes [24]. Berman et al. [27] introduced several models extending the FCLM to enable flows to deviate from predetermined shortest paths. Kim and Kuby [24] further adapted the FRLP to incorporate drivers' willingness to deviate from their shortest path for refueling, assuming unlimited station capacity. This enhanced model, known as the deviation-flow refueling location model, accounts for such deviations in its optimization framework.

Huang et al. [28], and Wang and Lin [29] addressed the problem of locating alternative fuel stations, assuming AFV users could utilize multiple deviation paths between all OD pairs in the network. They proposed a flow-based set covering model to minimize alternative fuel station location costs while satisfying travel demands between all OD pairs. Expanding on this, Upchurch et al. [30] extended the FRLP to account for the capacity limitations of refueling stations, resulting in the capacitated flow refueling location model. Furthermore, Hosseini and MirHassani [31] introduced a flexible formulation for the capacitated FRLP solving it with a heuristic method based on Lagrangian relaxation. In another study, Hosseini and MirHassani [32] considered the capacitated refueling station location problem with deviation from shortest paths. They proposed a mixed-integer linear programming model and a heuristic approach to solve the problem. Moreover, in another study, Hosseini and MirHassani [33] considered a stochastic variant of the capacitated FRLP, assuming uncertain traffic flow volumes on the network. They formulated the problem as a two-stage stochastic programming model, where the first stage involves locating permanent stations and the second stage involves placing portable stations. Similarly, Cheng and Wang [34] proposed a stochastic chance-constrained programming model for AFV, incorporating uncertain traffic flow. Their study addresses the challenges of maximizing alternative fuel usage and meeting travel demands in a traffic network by optimizing the locations of capacitated alternative fuel stations under budget constraints. Unlike, Hosseini and MirHassani [33] and Cheng and Wang [34], Mahmutoğullari and Yaman [35] addressed flow uncertainty by introducing a robust model for the refueling station location problem with routing (RSLP-R). This RSLP-R model is a maximal coverage approach that seeks to locate alternative fuel stations on a road network to maximize refueled AFV flows, considering both the limited range of vehicles and drivers' willingness to deviate for refueling.

Yildiz et al. [36] extended FRLP by incorporating the routing aspect of individual drivers which is known as deviation flow refueling location model-enhanced (DFRLM-E) [37]. They developed a mathematical model and solved it by a branch and price algorithm, which implicitly accounts for deviation tolerances without route pre-generation. Arslan et al. [38] tackled the refueling station location problem with routing (RSLP-R), aiming to maximize the flow of AFVs on a road network by determining optimal alternative refueling station locations while considering vehicle range limitations and driver deviation tolerances. Nordlund et al. [39] proposed an extended version of RSLP-R, known as the capacitated refueling station location problem with routing (CRSLP-R), which includes station capacities to limit the number of vehicles that can be refueled at each station. To address the complexity of CRSLP-R, they introduced two optimization methods: an expanded branch-and-cut approach and a branch-cut-and-price algorithm utilizing variables associated with feasible routes. Purba et al. [40] proposed hop constraints for the location-routing problem. These constraints limit the number of links a vehicle can travel before requiring a refueling stop. They aim to optimize the placement of emergency refueling stations for each alternative fuel type to facilitate evacuation routing. Their study introduces a Benders-inspired decomposition method utilizing a transformed network to address the location problem, and a metaheuristic branch-and-price approach for designing evacuation routes.

Chung and Kwon [41] proposed a model that addresses the gradual deployment of charging infrastructure over time, known as the multi-period flow-refueling location model (M-FRLM). Miralinaghi et al. [42] extended the M-FRLM approach by incorporating multi-period demand patterns and representing marginal operational costs through a staircase function. Their model focuses on locating capacitated alternative refueling stations to minimize the costs of construction and customer travel. To achieve this, they developed a mathematical model and applied branch-and-bound and Lagrangian relaxation algorithms. This extension not only aligns with the multi-period framework of M-FRLM but also addresses capacity constraints, enhancing the applicability of multi-period planning in infrastructure development. Additionally, Li et al. [43] demonstrated that multi-path and multi-period optimization principles could effectively determine the optimal locations for public electric charging stations, introducing this approach as the multi-path multi-period refueling location model. Medina et al. [44] developed an optimization model for EV fleet charging location assignment that addresses infrastructure planning challenges under electric grid constraints. Their model integrates substation capacity limits, operational characteristics of commercial fleets, and a flexible incentive-based mechanism to ensure grid stability and effective vehicle allocation.

Due to the limitations of discrete approaches in guaranteeing optimality for complex location problems, Abbaas and Ventura [8] used the continuous approach to solve the deviation location problem on a general network. They developed a polynomial-time algorithm to create a FDS and locate a single refueling station with the objective of maximizing flow coverage. Later they used network decomposition to break the problem into a set of easier sub problems improving the algorithm efficiency. In another study, Abbaas and Ventura [26] used lexicographic optimization to maximize coverage and minimize travel distance.

### *2.2. Research Gaps*

Most of the reviewed literature considers a predetermined finite set of candidate locations to establish the refueling stations, this is known as the discrete approach. Only a few studies consider an infinite set of candidate locations assuming that refueling stations can be located anywhere in the network, this is known as the continuous approach. A major limitation of the discrete approach is the lack of a theoretical guarantee for global optimality. The finite set is usually established to reduce the computational effort without considering solution qaltiy. Therefore, an optimal solution found using the discrete approach can be sub-optimal under the continuous approach [45]. To address this limitation, Kuby and Lim [46] and Kuby et al. [47] proposed three methods for generating additional candidate locations to better approximate the continuous solution space. Their mid-path segments method identifies arc segments where a single refueling station could effectively serve an entire path, reducing the need for multiple stations at vertices. Additionally, their Added Node Dispersion Problem (ANDP) uses heuristic dispersion-based strategies to distribute candidate locations along arcs, using maximin and minimax techniques. While these heuristic methods improve station placement compared to static node-based methods, they do not explicitly construct a FDS and, thus, lack a theoretical guarantee that a global optimal solution resides within the generated candidate set. Moreover, ANDP appraoch is computationally expensive. To improve the computational complexity, an extended version of the FRLP was proposed by Ventura et al. [48] and Kweon et al. [45]. They developed a polynomial-time algorithm for the continuous location problem on a tree network. The objective was to maximize the covered flow in daily round trips by a single refueling station. In these two papers, optimal locations were defined in terms of endpoints of refueling segments instead of taking midpoints of arcs as suggested by Kuby and Lim [46]. Later, Abbaas and Ventura [8, 26] considered the problem on a general network and allowed deviations from preplanned paths. Their study considers the characteristics of the traffic network and vehicle driving range to discretize the continuous version of the problem and find a FDS of candidate locations that guarantees optimality. However, they do not consider fleets with dedicated routes and stops along the way, leaving a clear research gap and a practical application unaddressed.

In this study, we consider vehicles with dedicated routes and service stops in a general network. These routes could correspond to public transportation vehicles, product deliveries from distribution centers to retail stores, or mail and package delivery. Vehicles are allowed to deviate from their preplanned routes for refueling services; however, deviation distance is limited. We leverage the problem structure and constraints to derive necessary and sufficient conditions for a candidate location to cover the flow on a dedicated route. The derived coverage conditions are then used to develop an exact algorithm to extract a FDS. We prove that at least one optimal solution exists within this set. Then we propose a mathematical model to minimize the number of refueling stations required to cover all routes in the network. Therefore, our methodology combines the computational efficiency of the discrete approach by creating a finite set of candidate locations, with the guaranteed optimality of the continuous approach by ensuring the dominance of our set of candidate locations.

## 3. Problem Statement

Consider a simple, connected, and undirected general network, $G(V, E)$, where $V$ is a set of $n$ vertices corresponding to network interchanges, and $E$ is a set of $e$ edges corresponding to road segments. An edge connecting vertices $v_i \in V$ and $v_j \in V$ is denoted by $(v_i, v_j) \in E$. A point in the network is defined as a location along an edge in $E$ or a vertex in $V$. A closed line segment between any two points $x_1, x_2$ on the same edge is denoted by $C(x_1, x_2)$, and the interior points of this line segment are referred to by int $(C(x_1, x_2))$. In a general network, multiple paths may exist between any two points $x_1, x_2$. The shortest path between these points, denoted by $P(x_1, x_2)$, is defined as the sequence of line segments and edges connecting them. The length of this path is denoted by $l(x_1, x_2)$.

In this study we propose a continuous approach algorithm to locate one or more refueling stations in a network with dedicated routes. These routes could be used by public transportation fleets, mail services, distribution operations, etc., to visit a sequence of stops in a certain order. These stops could be bus stops or delivery points, and they may be visited by multiple routes in the network. In this study, vehicles are allowed to deviate from their routes for refueling services. This allows a refueling station to serve multiple routes that do not share any common edges or vertices which reduces the total cost and improves the utilization of a refueling station. However, it may cause inconveniences due to the added travel distance to use a refueling station. Therefore, the additional driving distance is limited by a maximum value to maintain fuel efficiency and minimize inconvenience. The objective is to minimize the number of refueling stations that can cover all routes in the network.

Let $H = \{U_1, U_2, \dots, U_h\}$ be the set of routes in the network, where each route, $U_t \in H$, $t \in \{1,2,\dots,h\}$, is a directed closed walk that consists of a finite sequence of edges. Let $V(U_t)$ denote the set of vertices incident to the edges in $U_t$. The number of vehicles dedicated to route $U_t$ is referred to as the flow of route $U_t$ and is denoted by $f(U_t)$. Let $S_{U_t}$ be the finite sequence of stops along route $U_t$, where each stop $s_p \in S_{U_t}$ corresponds to either a vertex in $V(U_t)$ or a point along an edge in $U_t$. An edge may appear multiple times in a route in both directions. Therefore, an edge $(v_i, v_j) \in U_t$ will be written as $(v_i, v_j)$ to indicate travel direction from $v_i$ to $v_j$, or $(v_j, v_i)$ to indicate travel direction from $v_j$ to $v_i$. An intracity example with two routes is shown in Figure 2. These routes are $U_1 = \langle (v_2, v_3), (v_3, v_8), (v_8, v_{12}), (v_{12}, v_{16}), (v_{16}, v_{17}), (v_{17}, v_{16}), (v_{16}, v_{15}), (v_{15}, v_{14}), (v_{14}, v_{11}), (v_{11}, v_2) \rangle$, and $U_2 = \langle (v_5, v_8), (v_8, v_{12}), (v_{12}, v_9), (v_9, v_{10}), (v_{10}, v_6), (v_6, v_5) \rangle$. The length of $U_t$ is denoted by $l(U_t)$ and can be calculated as follows:

$$l(U_t) = \sum_{(v_i, v_j) \in U_t} l(v_i, v_j). \tag{1}$$

A sub-route $U_t(x_s, x_k)$ is defined as the part of route $U_t$ that starts at point $x_s$, located on edge $(v_i, v_s) \in U_t$ or at vertex $v_s$, and ends at point $x_k$, located on edge $(v_k, v_j) \in U_t$ or at vertex $v_k$, following the travel direction of $U_t$. Here, $v_s$ and $v_k$ are the first and last vertices encountered along sub-route $U_t(x_s, x_k)$, respectively. Therefore, $U_t(x_s, x_k)$ is a sequence that begins with the line segment $C(x_s, v_s)$, followed by a sequence of zero or more edges in $U_t$, and ends with the line segment $C(v_k, x_k)$. The length of sub-route $U_t(x_s, x_k)$ is denoted by $l(U_t(x_s, x_k))$ and can be calculated as follows:

$$l(U_t(x_s, x_k)) = l(x_s, v_s) + l(v_k, x_k) + \sum_{(v_i, v_j) \in U_t(v_s, v_k)} l(v_i, v_j). \tag{2}$$

Also, let $V(U_t(x_s, x_k))$ be the set of incident vertices of the edges in sub-route $U_t(x_s, x_k)$.

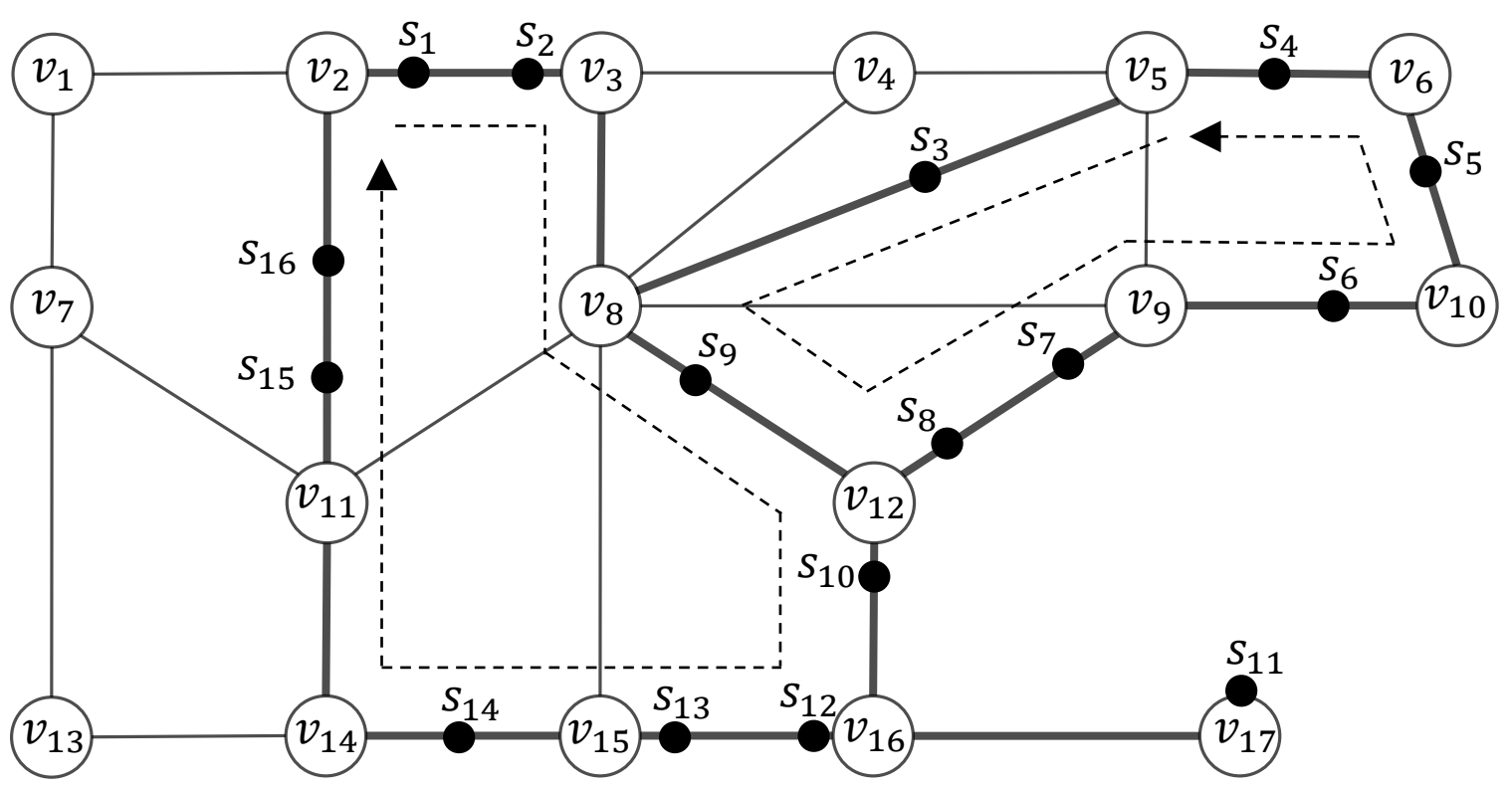

Figure 2: Intracity network

The following assumptions are used to formulate and solve the problem, these assumptions are consistent with the related literature [7, 26, 46, 49].

1. Refueling stations have unlimited capacities and are accessible from both sides of the road. This assumption ensures that all vehicles reaching a refueling station from any side of the road get refueling services.
2. All vehicles have the same driving range, denoted by $R$, for a full fuel tank. Fuel consumption is linearly correlated with driving distance. This assumption ensures that all vehicles can travel the same distance per refueling. Road conditions, congestions, slopes, and all other variables do not affect the driving range per full fuel tank.
3. Edges in network $G(V, E)$ are undirected, allowing vehicles to travel in both directions along any edge. As a result, the distance between any two points is identical in both directions, i.e., $l(x_1, x_2) = l(x_2, x_1)$. This assumption simplifies the discussion in this study. Additionally, although it is relatively common to have one-way roads within cities, we expect the lengths of edges in the network to be much shorter than $R$ in most intracity applications. Therefore, variations between $l(x_1, x_2)$ and $l(x_2, x_1)$ are assumed to be negligible in practice. This is consistent with the findings of Gayah and Daganzo [49].
4. The length of any route is at most $R$, vehicles can complete at least one trip around their dedicated route before running out of fuel. This assumption ensures that any route can be covered by a single refueling station located along one of its edges or their vertices.
5. Vehicles can deviate from their dedicated routes for refueling.
6. Vehicles are allowed to deviate at most a distance $D$ from their route to a refueling station location. The maximum allowed deviation distance can be defined as a fixed value, a percentage of the driving range $R$, or a percentage of the route length. Here, we will assume that $D$ is a fixed value and it is the same for all routes. Although $D$ is an upper deviation limit, if the vehicle completed one or more trips around its route but making an additional trip would cause the fuel level to drop below the level necessary to deviate a distance $D$ from the route, the driver must refuel as soon as possible. Figure 3 shows an example of this case; here, the length of the route is 14 distance units while $R = 100$. Starting with a full fuel tank, a vehicle can complete six trips around the route, but completing a seventh trip will only leave fuel enough to drive 2 additional distance units which is less than $D = 4$. Therefore, the vehicle should refuel after the sixth trip. Note that this assumption does not eliminate the case where $l(U_t) \leq R < l(U_t) + D$. In this case the vehicle cannot deviate the full distance $D$ after its first trip due to insufficient fuel.

7. U-turns are not allowed inside edges except at refueling stations. This assumption prohibits solutions that suggest making a U-turn in the middle of the road since this is impractical for most intracity applications. However, if there is a refueling station, we assume that vehicles can enter the station and then leave in any direction.
8. Triangle inequality is satisfied. This assumption ensures that the shortest path between any two points $x_1, x_2$ on the same edge $(v_i, v_j) \in E$ goes through the edge itself.

Based on these assumptions, let us define the set of routes that need to be considered to find an optimal solution. Any route $U_t$ where $l(U_t) > R$ cannot be covered by a single refueling station, this case is unlikely to happen in an intracity application. Also, if $f(U_t) = 0$ then $U_t$ does not affect the optimal solution. Therefore, the set of routes that needs to be examined for coverage is defined as follows:

$$H = \{U_t | l(U_t) \leq R, f(U_t) > 0\}. \quad (3)$$

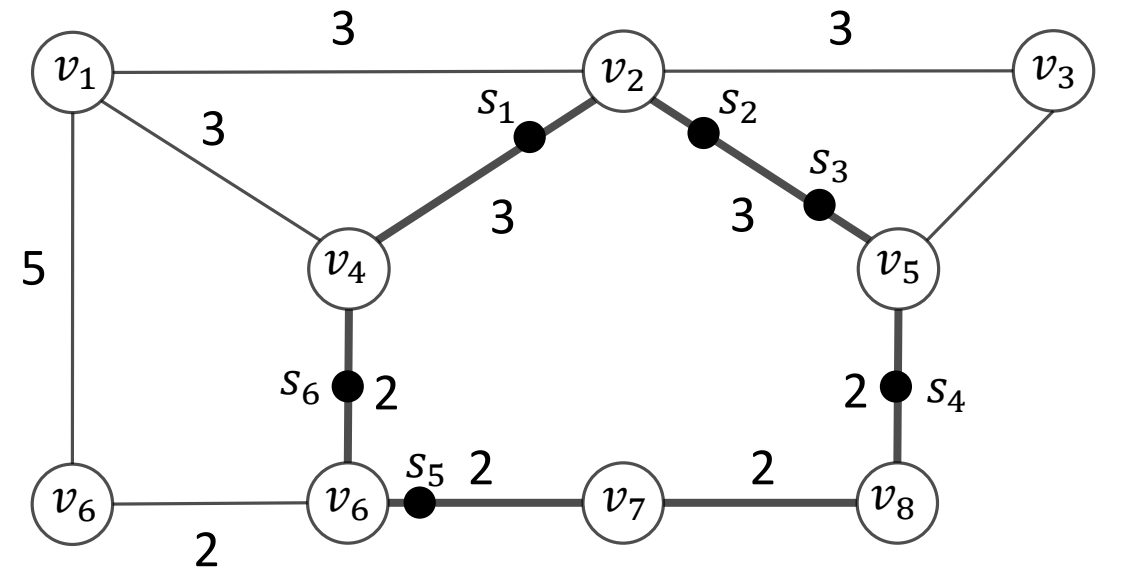

$U_1 = \langle (v_2, v_5), (v_5, v_8), (v_8, v_7), (v_7, v_6), (v_6, v_4), (v_4, v_2) \rangle$
$l(U_1) = 14, D = 4, R = 100$

Figure 3: Assumption 6 example

**4. Edge Scanning (ES) Algorithm for Dedicated Routes**

In this section, we introduce our proposed methodology to solve the continuous deviation flow refueling station location problem for dedicated routes. In Subsection 4.1, we derive the constraints that must be satisfied for a refueling station at a given candidate location to cover a certain flow, then we represent these constraints mathematically. In Subsection 4.2, we use the flow coverage constraints to identify all points and line segments in the network that can hold a refueling station capable of covering a given route/flow. These points and line segments are referred to as refueling sets. The following subsection presents an exact algorithm to find the refueling sets of all flows in the network. Lastly, we present a mathematical model to find the optimal combination of candidate locations that minimizes the number of required refueling stations to cover all flows in the network.

### *4.1. Flow Coverage Constraints*

Vehicles dedicated to a route, $U_t$, travel on the route's edges throughout the day until they need to refuel. Given that $l(U_t) \leq R$, if the refueling station is located along the edges of $U_t$ or their incident vertices, vehicles dedicated to route $U_t$ can refuel without deviating from their route. Therefore, any point $x$ along the edges of $U_t$ and their incident vertices can hold a refueling station that covers at least the flow $f(U_t)$. Since deviations are allowed, a refueling station serving $U_t$ may be located outside the route's edges and their incident vertices. In this case, let $s^l_{U_t}(x)$ be the last visited stop on route $U_t$ before deviating toward a refueling station located at point $x$. Also, let $s^f_{U_t}(s_p)$ be the next stop after stop $s_p$, following the defined travel direction of route $U_t$, $s_p \in S_{U_t}$. Additionally, let us define $v^a_{U_t}(s_p)$ as the vertex coinciding or immediately after stop $s_p$, and $v^b_{U_t}(s_p)$ as the vertex coinciding or immediately before stop $s_p$, following the defined travel direction of route $U_t$, $s_p \in S_{U_t}$. In order to minimize the total travel distance, we define $s^l_{U_t}(x)$ as follows:

$$s^l_{U_t}(x) = \arg\min_{s_p \in S_{U_t}} \{ l(s_p, v^a_{U_t}(s_p)) + l(v^a_{U_t}(s_p), x) + l(x, v^b_{U_t}(s^f_{U_t}(s_p))) + l(v^b_{U_t}(s^f_{U_t}(s_p)), s^f_{U_t}(s_p)) + l(U_t(s^f_{U_t}(s_p), s_p)) \}, \text{ ties are broken arbitrarily.} \quad (4)$$

Several cases can happen when deviating from a pre-planned route for refueling purposes. These cases affect when the deviation distance from the route begins to count against the maximum allowed deviation distance, as specified in assumption (6). Case 1 happens when the deviation path that minimizes the total travel distance leaves route $U_t$ using a subsequent vertex to $s^l_{U_t}(x)$ and goes toward point $x$, then returns to $U_t$ using a vertex preceding the next unvisited stop, $s^f_{U_t}(s^l_{U_t}(x))$. Figure 4 shows two examples illustrating Case 1. In both examples the refueling station is located at point $x$ on vertex $v_1$, $s^l_{U_t}(x) = s_5$, $s^f_{U_t}(s_5) = s_1$, and vehicles return to the route via $v^b_{U_t}(s_1) = v_2$. In Case 1.a, shown in Figure 4 part (a), vehicles deviate from the pre-planned route via $v_7$ which is the same as $v^a_{U_t}(s_5)$. However, in Case 1.b, shown in Figure 4 part (b), vehicles deviate from the pre-planned route via $v_8$ which is not the same as $v^a_{U_t}(s_5)$. Since the deviation from the pre-planned route in Case 1.a starts at $v^a_{U_t}(s_5) = v_7$, we start counting deviation distance against $D$ from $v_7$. However, in Case 1.b vehicles do not start deviating from their pre-planned route until they reach the deviation vertex $v_8$. Therefore, the distance from $v_7$ to $v_8$ is not counted against $D$.

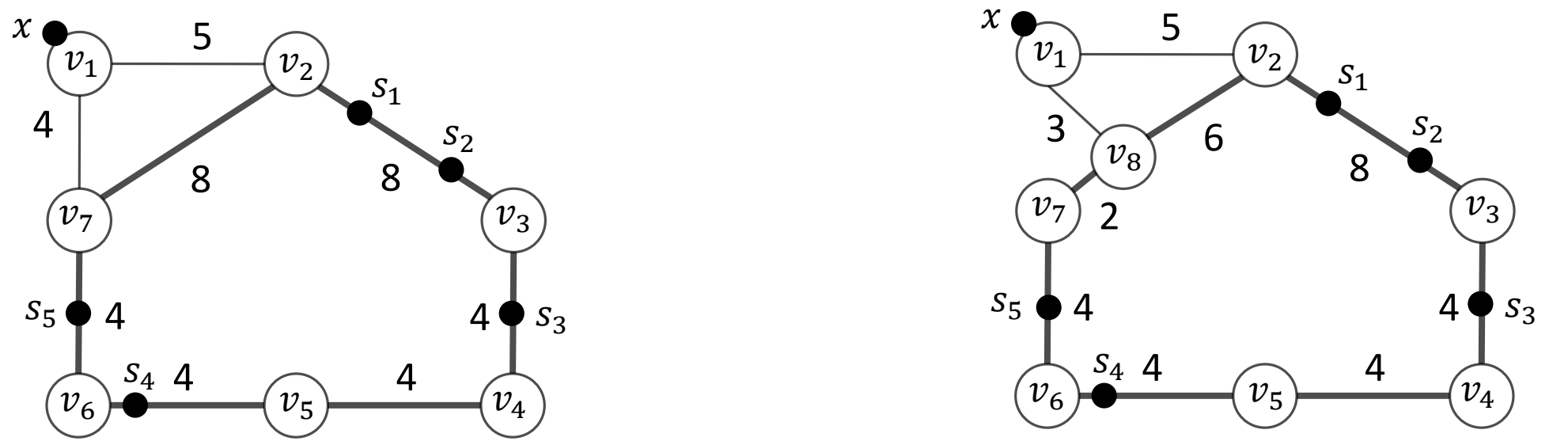

$U_1 = \langle(v_2, v_3), (v_3, v_4), (v_4, v_5), (v_5, v_6), (v_6, v_7), (v_7, v_2)\rangle$

(a) Case 1.a, the deviation vertex is $v_7 = v^a_{U_t}(s^l_{U_t}(x))$

$U_1 = \langle(v_2, v_3), (v_3, v_4), (v_4, v_5), (v_5, v_6), (v_6, v_7), (v_7, v_8), (v_8, v_2)\rangle$

(b) Case 1.b, the deviation vertex is $v_8 \neq v^a_{U_t}(s^l_{U_t}(x))$

Figure 4: Deviation Case 1

Case 2 happens when vehicles leave and/or return to their route using vertices that are not between the last visited stop $s^l_{U_t}(x)$ and the following stop $s^f_{U_t}(s^l_{U_t}(x))$. Figure 5 shows two examples illustrating Case 2. In both examples the refueling station is located at point $x$. Case 2.a, shown in Figure 5 part (a), happens when the deviation path that minimizes the total travel distance leaves $U_t$ using a vertex preceding $s^l_{U_t}(x)$, following the travel direction of route $U_t$. Case 2.b, shown in Figure 5 part (b), happens when the deviation path that minimizes the total travel distance returns using a subsequent vertex to $s^f_{U_t}(s^l_{U_t}(x))$, following the travel direction of route $U_t$.

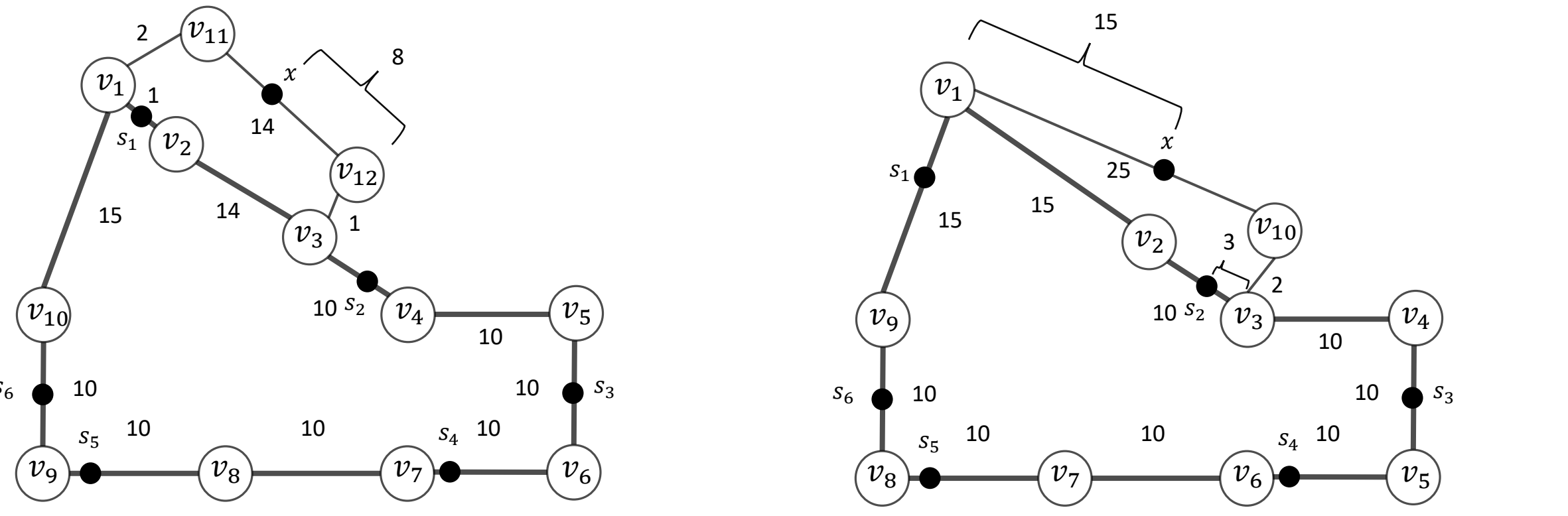

$U_1 = \langle(v_1, v_2), (v_2, v_3), (v_3, v_4), (v_4, v_5), (v_5, v_6), (v_6, v_7), (v_7, v_8), (v_8, v_9), (v_9, v_{10}), (v_9, v_{10})\rangle$

(a) Case 2.a deviation vertex precedes $s^l_{U_t}(x)$

$U_1 = \langle(v_1, v_2), (v_2, v_3), (v_3, v_4), (v_4, v_5), (v_5, v_6), (v_6, v_7), (v_7, v_8), (v_8, v_9), (v_9, v_1)\rangle$

(b) Case 2.b, return vertex is subsequent to $s^f_{U_t}\left(s^l_{U_t}(x)\right)$

Figure 5: Deviation Case 2

In Figure 5 part (a), $s^l_{U_t}(x) = s_1$, $v^a_{U_t}(s_1) = v_2$, $s^f_{U_t}(s_1) = s_2$, and $v^b_{U_t}(s_2) = v_3$. The deviation path that minimizes the total travel distance, including the refueling trip, does not go through edge $(v_2, v_3)$.

Rather, after visiting $s_1$ it goes back to $v_1$ then toward point $x$ via $v_{11}$. However, due to practical and legal considerations we assume that vehicles cannot make a U-turn inside an edge except at refueling locations, assumption 7. Therefore, vehicles have to reach a vertex before they can reverse their travel direction. Based on that, when a vehicle in Figure 5 part (a) needs to refuel, it visits $s_1$ then continues on edge $(v_1, v_2)$ until it reaches $v_2$. After that it makes a U-turn to go back to $v_1$, and then goes toward point $x$ via $v_{11}$. The vehicle returns to the route $U_1$ via vertices $v_{12}$ and $v_3$. In Figure 5 part (b), $s^l_{U_t}(x) = s_1$, $v^a_{U_t}(s_1) = v_1$, $s^f_{U_t}(s_1) = s_2$, and $v^b_{U_t}(s_2) = v_2$. Here, a vehicle visits $s_1$ then leaves the route using $v_1$ toward point $x$. After reaching point $x$, it returns to the route $U_1$ using $v_3$ and goes toward the next unvisited stop $s_2$. However, to resume following the defined travel direction of the route, the vehicle needs to make a U-turn. Therefore, it must reach vertex $v_2$ then make a U-turn.

To calculate the deviation distance traveled by vehicles, let us define $v^d_{U_t}(x)$ as the deviation vertex where vehicles start deviating from their route $U_t$, either by leaving its edges or by going against its defined travel direction. In Cases 1.a and 1.b by definition $v^d_{U_t}(x)$ must be the same as $v^a_{U_t}(s^l_{U_t}(x))$ or a subsequent vertex up to $v^b_{U_t}(s^f_{U_t}(s^l_{U_t}(x)))$. This can be seen in Figure 4 parts (a) and (b) where $v^d_{U_t}(x) = v^a_{U_t}(s^l_{U_t}(x)) = v_7$ and $v^d_{U_t}(x) = v_8$, respectively. Also, in Figure 5 part (b) where $v^d_{U_t}(x) = v^a_{U_t}(s^l_{U_t}(x)) = v_1$. In Case 2.a, $v^d_{U_t}(x)$ is the same as $v^a_{U_t}(s^l_{U_t}(x))$ since vehicles reach $v^a_{U_t}(s^l_{U_t}(x))$ and make a U-turn to go against the defined travel direction of the route. This can be seen in Figure 5 part (a). Based on that, the deviation vertex, $v^d_{U_t}(x)$, must be the same as $v^a_{U_t}(s^l_{U_t}(x))$ or a subsequent vertex up to $v^b_{U_t}(s^f_{U_t}(s^l_{U_t}(x)))$. Additionally, it must be an incident vertex to an edge that belongs to the shortest path toward the refueling station, where this edge does not belong to the sub-route between $v^a_{U_t}(s^l_{U_t}(x))$ and $v^b_{U_t}(s^f_{U_t}(s^l_{U_t}(x)))$. Mathematically, we can define the deviation vertex $v^d_{U_t}(x)$ as follows:

$$v^d_{U_t}(x) = \{v_i | (v_i, v_j) \in P(v^a_{U_t}(s^l_{U_t}(x)), x), v_i \in V\left(U_t\left(v^a_{U_t}(s^l_{U_t}(x)), v^b_{U_t}(s^f_{U_t}(s^l_{U_t}(x)))\right)\right), v_j \notin V\left(U_t\left(v^a_{U_t}(s^l_{U_t}(x)), v^b_{U_t}(s^f_{U_t}(s^l_{U_t}(x)))\right)\right), V\left(U_t\left(v^a_{U_t}(s^l_{U_t}(x)), v^b_{U_t}(s^f_{U_t}(s^l_{U_t}(x)))\right)\right) \subseteq V(U_t), \text{and } V(U_t), v_j \subseteq V\}. \quad (5)$$

Note that, vertices $v^a_{U_t}(s^l_{U_t}(x))$, $v^b_{U_t}(s^f_{U_t}(s^l_{U_t}(x)))$, and $v^d_{U_t}(x)$ could be the same vertex, two different vertices, or three different vertices. Additionally, note that in all cases vehicles must reach $v^b_{U_t}(s^f_{U_t}(s^l_{U_t}(x)))$ after a refueling trip.

A route $U_t$ is considered to be covered by a single refueling station located at a point $x$ if any vehicle dedicated to this route can refuel at point $x$, complete at least one trip visiting all stops in the route, then go back to point $x$ without running out of fuel or exceeding the maximum allowed deviation distance, $D$. Based on that and by assumption (4), any point located on the edges of a route or at their incident vertices can simply cover the route. However, if point $x$ is not located along any edge in $U_t$ or at an incident vertex, then the following two constraints must be satisfied for $x$ to cover $U_t$:

$$l\big(v_{U_t}^d(x), x\big) + l\big(x, v_{U_t}^b(s_{U_t}^f(s_{U_t}^l(x)))\big) + l'\big(U_t(v_{U_t}^b(s_{U_t}^f(s_{U_t}^l(x))), v_{U_t}^d(x))\big) \le R, \tag{6}$$

$$l\big(v_{U_t}^d(x), x\big) \le D. \tag{7}$$

where,

$$l'\big(U_t(v_{U_t}^b(s_{U_t}^f(s_{U_t}^l(x))), v_{U_t}^d(x))\big) = \begin{cases} l(U_t), v_{U_t}^b\left(s_{U_t}^f\left(s_{U_t}^l(x)\right)\right) = v_{U_t}^d(x), \\ l\big(U_t(v_{U_t}^b(s_{U_t}^f(s_{U_t}^l(x))), v_{U_t}^d(x))\big), \text{otherwise.} \end{cases} \tag{8}$$

Constraint (6) guarantees that the length of the entire trip from the refueling station and back to it is less than or equal to the vehicle's driving range, $R$. It ensures that vehicles can always come back to the refueling station without running out of fuel. Here we define $l'\big(U_t(v_{U_t}^b(s_{U_t}^f(s_{U_t}^l(x))), v_{U_t}^d(x))\big)$ to account for the case when $v_{U_t}^b(s_{U_t}^f(s_{U_t}^l(x)))$ and $v_{U_t}^d(x)$ are the same vertex. In this case, $l'\big(U_t(v_{U_t}^b(s_{U_t}^f(s_{U_t}^l(x))), v_{U_t}^d(x))\big)$ equals the length of the route, $l(U_t)$. Otherwise, $l'\big(U_t(v_{U_t}^b(s_{U_t}^f(s_{U_t}^l(x))), v_{U_t}^d(x))\big)$ equals the length of the sub-route from $v_{U_t}^b(s_{U_t}^f(s_{U_t}^l(x)))$ to $v_{U_t}^d(x))$ following the defined travel direction of route $U_t$, i.e. $l\big(U_t(v_{U_t}^b(s_{U_t}^f(s_{U_t}^l(x))), v_{U_t}^d(x))\big)$. Constraint (7) guarantees that the deviation distance to point $x$ cannot exceed $D$.

The first flow coverage constraint, (6), is necessary when the route is relatively long and the amount of fuel remaining after completing one full trip around the route may not be enough to deviate the full distance $D$. An example of this case is shown in Figure 6. Here, the driving range $R = 20$ while the length of the route is $l(U_1) = 18$. Note that, any point within 1 distance unit from vertex $v_8$ along edge $(v_8, v_9)$ satisfies the two constraints and can cover route $U_1$. However, points beyond that cannot cover the route since the total length of a complete trip will exceed $R$, even if these points are within distance $D$ from the deviation vertex $v_8$. For example, take point $y$, where $l(v_8, y) = 2$. The length of a complete trip starting from point $y$, then visiting all stops around the route, and coming back to point $y$ to refuel would be $l(y, v_8) + l(U_t) + l(v_8, y) = 22$. Therefore, a vehicle refueling at point $y$ will run out of fuel before reaching the refueling station again, even though $d(v_8, y) < D$. Additionally, even though the second flow coverage constraint, (7), only limits the deviation distance from the deviation vertex toward the refueling station location explicitly, it implicitly limits the total additional driving distance for refueling purposes to be less than or equal to $2D$. This result is presented in Lemma 1.

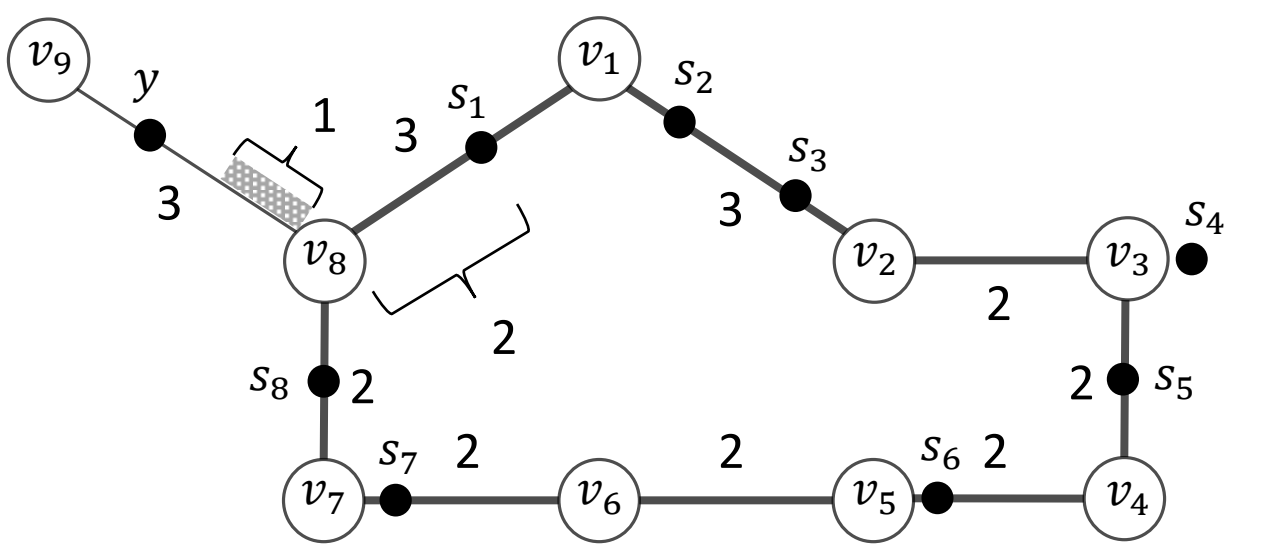

$U_1 = \langle (v_1, v_2), (v_2, v_3), (v_3, v_4), (v_4, v_5), (v_5, v_6), (v_6, v_7), (v_7, v_8), (v_8, v_1) \rangle$
$l(U_1) = 18, D = 4, R = 20$

Figure 6: Flow coverage constraint (6) is binding

**Lemma 1.** Let $U_t$ be a dedicated route, $v_{U_t}^d(x)$ is the deviation vertex to a refueling station located at point $x$ which does not lie on any edge in $U_t$ or their incident vertices. Given that the deviation distance from $v_{U_t}^d(x)$ to $x$ is bounded by $D$ and drivers take the shortest path from $x$ to the next unvisited stop on $U_t$, the total additional driving distance beyond $l(U_t)$ per refueling trip cannot exceed $2D$.

**Proof.** Consider a vehicle traveling on route $U_t$. After visiting stop $s_{U_t}^l(x)$, it deviates starting from vertex $v_{U_t}^d(x)$ to refuel at point $x$ and then goes to the next unvisited stop, $s_{U_t}^f(s_{U_t}^l(x))$, passing through vertex $v_{U_t}^b(s_{U_t}^f(s_{U_t}^l(x)))$. The return trip takes the shortest path from $x$ to $v_{U_t}^b(s_{U_t}^f(s_{U_t}^l(x)))$ which may or may not pass through the deviation vertex $v_{U_t}^d(x)$. Therefore, the length of this return trip must satisfy the following Inequality:

$$l(x, v_{U_t}^b(s_{U_t}^f(s_{U_t}^l(x)))) \leq l(x, v_{U_t}^d(x)) + l(U_t(v_{U_t}^d(x), v_{U_t}^b(s_{U_t}^f(s_{U_t}^l(x))))), \tag{9}$$

Given that $l(v_{U_t}^d(x), x) \leq D$, then the length of the return trip must be:

$$l(x, v_{U_t}^b(s_{U_t}^f(s_{U_t}^l(x)))) \leq D + l(U_t(v_{U_t}^d(x), v_{U_t}^b(s_{U_t}^f(s_{U_t}^l(x))))), \tag{10}$$

The total additional driving distance beyond $l(U_t)$, denoted by $AD$, can be calculated as follows:

$$AD = \Big(l\big(v_{U_t}^d(x), x\big) + l\big(x, v_{U_t}^b(s_{U_t}^f(s_{U_t}^l(x)))\big)\Big) - l\,\Big(U_t\big(v_{U_t}^d(x), v_{U_t}^b(s_{U_t}^f(s_{U_t}^l(x)))\big)\Big), \tag{11}$$

Substituting Inequality (10) into Equation (11) we get:

$$\begin{aligned} AD \leq l\big(v_{U_t}^d(x), x\big) + D + l(U_t(v_{U_t}^d(x), v_{U_t}^b(s_{U_t}^f(s_{U_t}^l(x))))) - \\ l(U_t(v_{U_t}^d(x), v_{U_t}^b(s_{U_t}^f(s_{U_t}^l(x))))), \end{aligned} \tag{12}$$

Substituting Inequality (7) into (12) we get:

$$AD \leq 2D. \square \tag{13}$$

Based on this discussion we can define the set of routes in $H$ covered by any point $x$ as follows:

$$T(x) = \{U_t | U_t \in H, [x \in C(v_i, v_j), (v_i, v_j) \in U_t] \; or \; [x \notin C(v_i, v_j) \; \forall \; (v_i, v_j) \in U_t, \; l(v_{U_t}^d(x), x) + l(x, v_{U_t}^b(s_{U_t}^f(s_{U_t}^l(x)))) + l'(U_t(v_{U_t}^b(s_{U_t}^f(s_{U_t}^l(x))), v_{U_t}^d(x))) \leq R, l(v_{U_t}^d(x), x) \leq D]\}. \tag{14}$$

And the total flow covered by point $x$, denoted by $F(x)$, can be found using:

$$F(x) = \sum_{U_t \in T(x)} f(U_t). \tag{15}$$

*4.2. Refueling Sets*

A refueling set for route $U_t$ on edge $(v_a, v_b) \in E$, denoted by $RS(U_t; v_a, v_b)$, is the set of all points on edge $(v_a, v_b)$ that can cover route $U_t$. Therefore, $U_t \in T(x)$ for all $x \in RS(U_t; v_a, v_b)$. There are two cases to be considered when looking for refueling sets that cover a dedicated route $U_t \in H$ on an edge $(v_a, v_b) \in E$. Case (1) occurs when edge $(v_a, v_b)$ belongs to the route $U_t$, $(v_a, v_b) \in U_t$. Given the assumption that the length of any route is at most $R$, assumption (4), any point on edge $(v_a, v_b)$ can cover $U_t$. The refueling set for route $U_t$ on edge $(v_a, v_b) \in U_t$ is defined as follows:

$$RS(U_t; v_a, v_b) = \{x | x \in C(v_a, v_b), (v_a, v_b) \in U_t\}. \tag{16}$$

Case (2) occurs when the edge $(v_a, v_b) \in E$ does not belong to the route $U_t$, $(v_a, v_b) \notin U_t$. In this case, a vehicle needs to deviate from its dedicated route to reach a refueling station location on edge $(v_a, v_b)$. The vehicle may be able to enter the edge using any of its two vertices or only one of them. Therefore, to define the refueling set on edge $(v_a, v_b) \notin U_t$ we need to determine the remaining possible travel distance according to the flow coverage constraints (6) and (7) at the two possible entry vertices, $v_a$ and $v_b$. Let $\beta^{(q)}(U_t; v_a, v_b)$, $q \in \{a, b\}$, be the remaining possible travel distance along edge $(v_a, v_b)$ when vehicles of route $U_t$ access this edge via vertex $v_q$, $q \in \{a, b\}$, considering the first coverage constraint, (6), only. $\beta^{(q)}(U_t; v_a, v_b)$ can be defined as follows:

$$\beta^{(q)}(U_t; v_a, v_b) = \begin{cases} l(v_a, v_b), R - \begin{pmatrix} l(v_{U_t}^d(v_q), v_q) + l(v_w, v_{U_t}^b(s_{U_t}^f(s_{U_t}^l(v_q)))) + \\ l'(U_t(v_{U_t}^b(s_{U_t}^f(s_{U_t}^l(v_q))), v_{U_t}^d(v_q))) \end{pmatrix} \geq l(v_a, v_b), \\ \left( R - \begin{pmatrix} l(v_{U_t}^d(v_q), v_q) + l(v_q, v_{U_t}^b(s_{U_t}^f(s_{U_t}^l(v_q)))) + \\ l'(U_t(v_{U_t}^b(s_{U_t}^f(s_{U_t}^l(v_q))), v_{U_t}^d(v_q))) \end{pmatrix} \right) / 2, \text{Otherwise}, \end{cases} \tag{17}$$

$$q \in \{a, b\}; w \in \{a, b\}; w \neq q.$$

In Equation (17), the upper part represents the case when a vehicle can enter edge $(v_a, v_b) \notin U_t$ using one vertex $v_q$, $q \in \{a, b\}$ and reach the other vertex $v_w$, $w \in \{a, b\}$, $w \neq q$, then go back to the route to complete a full trip without running out of fuel. In this case, the entire edge satisfies the first constraint, (6). On the other hand, the lower part of Equation (17) represents the case when a vehicle enters edge $(v_a, v_b) \notin U_t$ using vertex $v_q$, $q \in \{a, b\}$, but cannot reach the other vertex $v_w$, $w \in \{a, b\}$, $q \neq w$, due to the limited driving range. Therefore, it has to make a U-turn at a potential refueling station and leave the edge using the same entry vertex $v_q$, $q \in \{a, b\}$. Note that, the distance $R - l(v_{U_t}^d(v_q), v_q) + l(v_q, v_{U_t}^b(s_{U_t}^f(s_{U_t}^l(v_q)))) + l'(U_t(v_{U_t}^b(s_{U_t}^f(s_{U_t}^l(v_q))), v_{U_t}^d(v_q)))$ is divided by 2 because vehicles need to travel from the entry vertex, $v_q$, to the refueling station location and then go back to $v_q$ to exit the edge. Figure 6 can be used to explain the calculation of $\beta^{(q)}(U_t; v_a, v_b)$. To find $\beta^{(8)}(U_t; v_8, v_9)$, we subtract the following distances from $R$: first, the distance from the deviation vertex, $v_{U_1}^d(v_8) = v_8$, to the entry vertex of edge $(v_8, v_9)$ which is $v_8$, that is $l(v_8, v_8) = 0$; second, the distance from the other vertex, $v_9$, to the vertex before the next unvisited stop $v_{U_1}^b(s_1)$ which is $v_8$, that is $l(v_9, v_8) = 3$; third, the length of the route, $l(U_1) = 18$. The result is $20 - 21 = -1 < l(v_8, v_9)$. This means that vehicles cannot go through the entire edge $(v_8, v_9)$, return to the route, and then complete one trip without running out of fuel. Therefore, we need to examine the second part of Equation (17). In the second part, subtract the following distances from $R$: first, the distance from the deviation vertex, $v_8$, to the entry vertex of edge $(v_8, v_9)$, that is $l(v_8, v_8) = 0$; second, the distance from the same entry vertex of the edge to the vertex before the next unvisited stop $v_{U_1}^b(s_1)$ which is $v_8$, that is $l(v_8, v_8) = 0$, and the length of the route, $l'(U_t(v_{U_t}^b(s_{U_t}^f(s_{U_t}^l(v_q))), v_{U_t}^d(v_q))) = l(U_1) = 18$. The result is $20 - 18 = 2$. The value of $\beta^{(8)}(U_1; v_8, v_9) = 2/2 = 1$. Therefore, the remaining possible travel distance according to the first flow coverage constraint, Constraint (6), when entering edge $(v_8, v_9)$ using vertex $v_8$ is one distance unit. Note that, $\beta^{(9)}(U_1; v_8, v_9) = -2$. A negative value for $\beta^{(q)}(U_t; v_a, v_b)$ means that vehicles attempting to complete a full trip around $U_t$ and then enter edge $(v_a, v_b)$ using $v_q$, $q \in \{a, b\}$, will run out of fuel before reaching $v_q$. $\beta^{(q)}(U_t; v_a, v_b) = 0$ means that vehicles attempting to complete a full trip around $U_t$ and then enter edge $(v_a, v_b)$ using $v_q$, $q \in \{a, b\}$, will not be able to reach any interior point of the edge beyond $v_q$ without running out of fuel.

Next, let $\delta^{(q)}(U_t; v_a, v_b)$, $q \in \{a, b\}$, be the remaining allowed deviation distance along edge $(v_a, v_b)$ when vehicles access this edge via vertex $v_q$, $q \in \{a, b\}$, considering the two flow coverage constraints, Constraints (6) and (7). $\delta^{(q)}(U_t; v_a, v_b)$ is defined as follows:

$$\delta^{(q)}(U_t; v_a, v_b) = \min\{D - l(v^d_{U_t}(v_q), v_q), \beta^{(q)}(U_t; v_a, v_b)\},\ q \in \{a, b\}. \tag{18}$$

Back to the example in Figure 6, given that $\beta^{(8)}(U_1; v_8, v_9) = 1$ and $D = 4$, then $\delta^{(8)}(U_1; v_8, v_9) = \min\{4 - l(v_8, v_8), 1\} = \min\{4,1\} = 1$. Note that, considering $v_9$ as the entry vertex for edge $(v_8, v_9)$ in this example gives $\delta^{(9)}(U_1; v_8, v_9) = \min\{4 - 3, -2\} = -2$. A negative value for $\delta^{(q)}(U_t; v_a, v_b)$ means that none of the points on edge $(v_a, v_b)$ satisfies the two coverage constraints when entering the edge using vertex $v_q$, $q \in \{a, b\}$. Figure *7* shows another example where constraint (7) is the binding constraint. Here, it is desired to find the refueling set that can cover route $U_1$ on edge $(v_a, v_b)$. Since $R$ is much longer than the length of the route, $\beta^{(a)}(U_t; v_a, v_b) = \beta^{(b)}(U_t; v_a, v_b) = l(v_a, v_b) = 5$. Next, if vertex $v_a$ is used to enter edge $(v_a, v_b)$ the deviation vertex from route $U_1$ would be $v_5$ and the remaining allowed deviation distance at $v_a$ would be 2 units, $\delta^{(a)}(U_1; v_a, v_b) = \min\{4 - l(v_5, v_a), 5\} = \min\{2,5\} = 2$. On the other hand, if the vehicle enters edge $(v_a, v_b)$ via vertex $v_b$ the deviation vertex from route $U_1$ would be $v_6$ and the remaining allowed deviation distance at $v_b$ would be 1 unit, $\delta^{(b)}(U_1; v_a, v_b) = \min\{4 - l(v_6, v_b), 5\} = \min\{1,5\} = 1$.

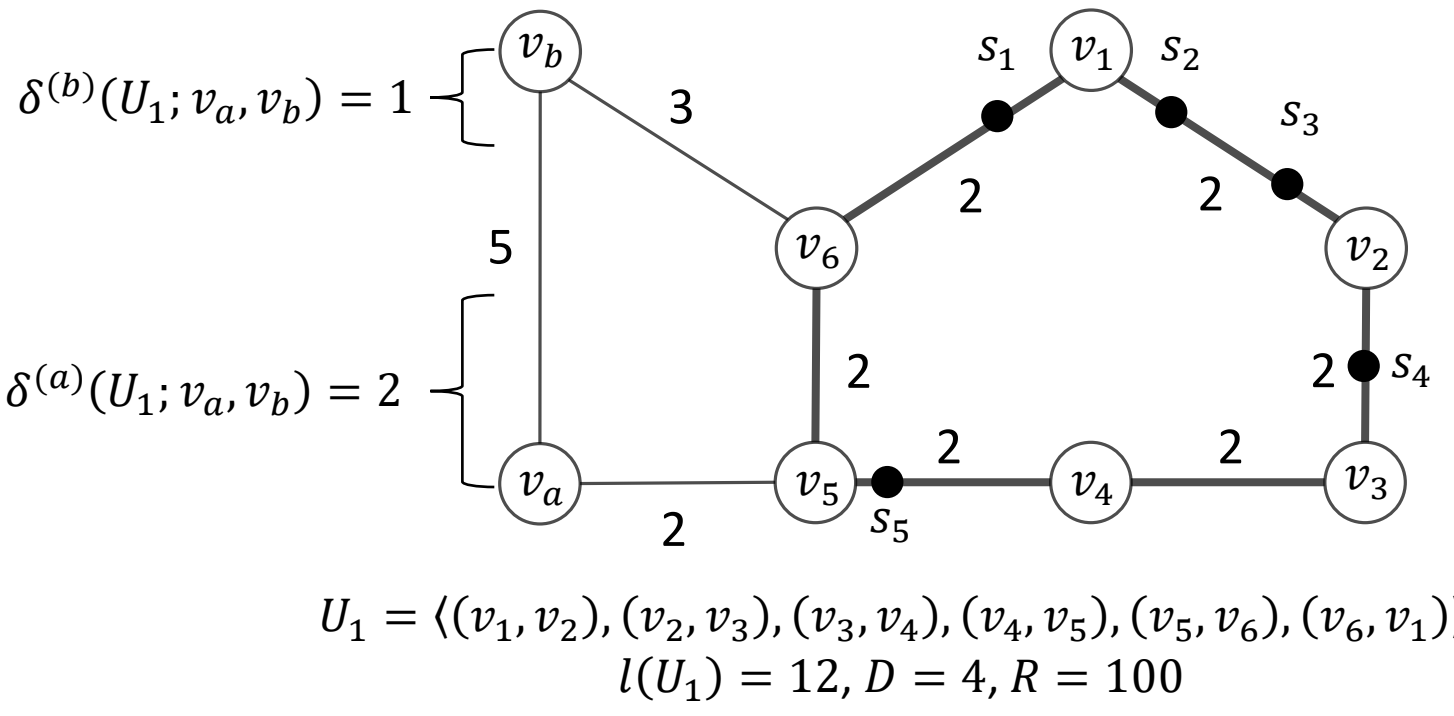

Figure 7: Refueling set

As can be seen in Figure 7 it is possible to have two disconnected line segments on edge $(v_a, v_b)$ that satisfy constraints (6) and (7); each line segment is defined assuming a different entry vertex. These two segments are called refueling segments and they are denoted by $RS^{(q)}(U_t; v_a, v_b)$, $q \in \{a, b\}$, depending on the entry vertex used in the definition of the segment.

$$RS^{(q)}(U_t; v_a, v_b) = \{x \mid x \in C(v_a, v_b), (v_a, v_b) \notin U_t, l(v_q, x) \le \delta^{(q)}(U_t; v_a, v_b)\},\ q \in \{a, b\}. \tag{19}$$

Now, to establish the overall refueling set that can cover route $U_t$ on edge $(v_a, v_b) \notin U_t$ the following equation is used:

$$RS(U_t; v_a, v_b) = RS^{(a)}(U_t; v_a, v_b) \cup RS^{(b)}(U_t; v_a, v_b). \tag{20}$$

Note that, if any of the two remaining deviation distances, $\delta^{(a)}(U_t; v_a, v_b)$ or $\delta^{(b)}(U_t; v_a, v_b)$, exceeds the length of the edge, i.e. if any of the two refueling segments $RS^{(a)}(U_t; v_a, v_b)$ or $RS^{(b)}(U_t; v_a, v_b)$ include the entire edge $(v_a, v_b)$, then there is no need to find the other one since $RS(U_t; v_a, v_b)$ will include all points in $C(v_a, v_b)$.

Each refueling segment $RS^{(q)}(U_t; v_a, v_b)$, $q \in \{a, b\}$ can have up to two endpoints depending on the value of $\delta^{(q)}(U_t; v_a, v_b)$. If $\delta^{(q)}(U_t; v_a, v_b)$ is negative, then vehicles dedicated to route $U_t$ cannot reach edge $(v_a, v_b)$ via vertex $v_q$. However, if $\delta^{(q)}(U_t; v_a, v_b)$ is more than or equal to zero, route $U_t$ has a non-empty refueling segment $RS^{(q)}(U_t; v_a, v_b)$ on edge $(v_a, v_b)$. If $\delta^{(q)}(U_t; v_a, v_b) = 0$ the refueling segment consists of a single point, that is vertex $v_q$; however, if $\delta^{(q)}(U_t; v_a, v_b) > 0$, the segment has two endpoints, one of them is $v_q$. Therefore, for any edge $(v_a, v_b) \in E$ the refueling set $RS(U_t; v_a, v_b)$ can have up to four endpoints, denoted by $w^k_{U_t;v_a,v_b}$, $k = 1, \dots, 4$. Recall that for the case where $(v_a, v_b) \in U_t$ any point on the edge belongs to $RS(U_t; v_a, v_b)$ and there are two endpoints, $v_a$ and $v_b$. In general, the set of endpoints of $RS(U_t; v_a, v_b)$ is denoted by $EP(U_t; v_a, v_b)$ and is defined as follows:

$$EP(U_t; v_a, v_b) = \{w^k_{U_t;v_a,v_b} | w^k_{U_t;v_a,v_b}, k = 1, \dots, 4, \text{ are endpoints of } RS(U_t; v_a, v_b)\}. \quad (21)$$

The set of endpoints on edge $(v_a, v_b) \in E$ and the overall set of endpoints in the network are denoted by $EP(v_a, v_b)$ and $EP$, respectively. These sets are defined as follows:

$$EP(v_a, v_b) = \cup_{U_t \in H} EP(U_t; v_a, v_b), \quad (22)$$

$$EP = \cup_{(v_a, v_b) \in E} EP(v_a, v_b). \quad (23)$$

Refueling sets of various routes may intersect on network edges as shown in Figure 8. In this case, the intersection set covers all of the routes that can be covered by the individual intersecting sets. The example in Figure 8 shows two routes $U_1$ and $U_2$, the maximum allowed deviation distance $D = 5$ units, and the length of each edge is shown in the figure. Here, $\delta^{(a)}(U_1; v_a, v_b) = 4$ and $\delta^{(b)}(U_1; v_a, v_b) = -1$; therefore, vehicles dedicated to route $U_1$ can reach edge $(v_a, v_b)$ only via vertex $v_a$ and the refueling set $RS(U_1; v_a, v_b) = \{x | x \in (v_a, v_b), l(v_a, x) \le 4\}$. Similarly, $\delta^{(a)}(U_2; v_a, v_b) = -3$ and $\delta^{(b)}(U_2; v_a, v_b) = 2$; therefore, vehicles dedicated to route $U_2$ can reach edge $(v_a, v_b)$ only via vertex $v_b$ and the refueling set $RS(U_2; v_a, v_b) = \{x | x \in (v_a, v_b), l(v_b, x) \le 2\}$. These refueling sets intersect in the shaded line segment on edge $(v_a, v_b)$; therefore, this segment satisfies the flow coverage constraints for the two routes and can cover both of them. Let the set of routes covered by endpoint $w \in EP$ be $T(w)$ and the total flow covered by $w$ be $F(w)$. $T(w)$ and $F(w)$ can be defined as follows:

$$T(w) = \{U_t | U_t \in H, w \in R(U_t; v_a, v_b), (v_a, v_b) \in E\}, \quad (24)$$

$$F(w) = \sum_{U_t \in T(w)} f(U_t). \tag{25}$$

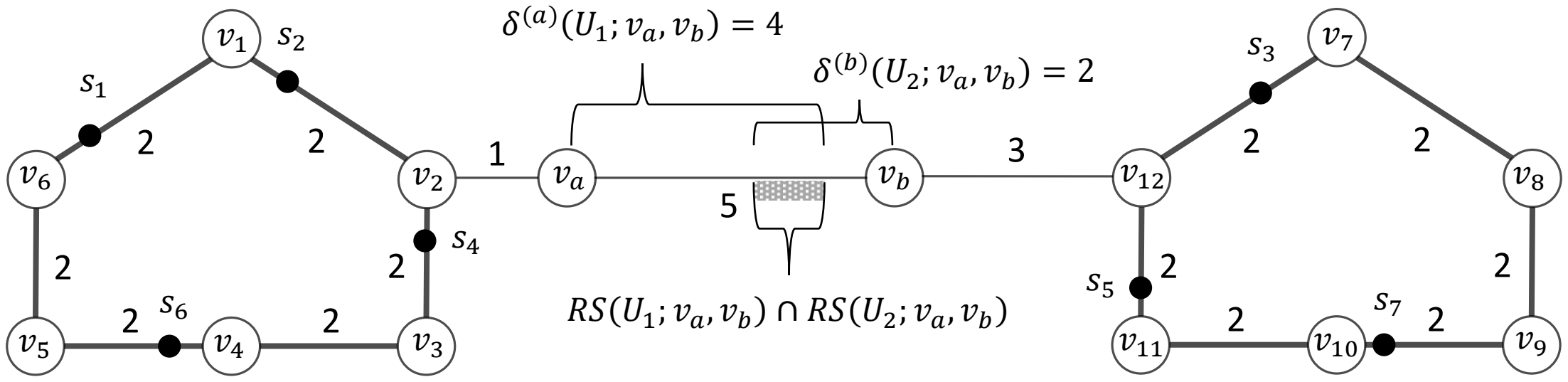

Figure 8: Intersecting refueling sets

Here, it is important to note that although a non-empty refueling set for any route on any edge must include at least one vertex, it is not sufficient to only use network vertices as candidate locations to find an optimal solution. In fact, many researchers studied similar problems starting with a predetermined set of candidate locations. These locations are usually network vertices, or a set that consists of both network vertices and some points along the edges like midpoints. These methods can lead to sub-optimal solutions. For example, if only network vertices were considered in Figure 8, at least two refueling stations will be required to cover the two routes. Even with a candidate location added in the middle of edge $(v_a, v_b)$ the minimum required number of stations to cover the two routes would still be two. However, with the continuous approach, only one refueling station is required in the intersection between $RS(U_1; v_a, v_b)$ and $RS(U_2; v_a, v_b)$.

Now, we extend Theorem 1 from [8] to show that the set of endpoints, $EP$, has at least one optimal solution to the continuous deviation flow refueling station location problem for dedicated routes.

**Theorem 1.** For any point $x$ located either at a vertex in $V$ or along an edge in $E$, there is always an endpoint $w \in EP$ such that $T(x) \subseteq T(w)$. Moreover, in the dedicated routes refueling station location problem, there is always an optimal solution $X^*$, where $X^* \subseteq EP$.

**Proof.** First, let $x$ be a point that does not belong to any refueling set. By definition, $x$ cannot cover any route, i.e., $T(x) = \emptyset$. This implies that $T(x) \subseteq T(w)$ for all $w \in EP$.

Next, let $x$ be a point that belongs to a refueling set along edge $(v_a, v_b) \in E$. In this case, $x$ could be an endpoint or an interior point. If $x$ is an endpoint, then the result follows. Otherwise, $x$ is an interior point. Note that, since refueling sets are defined within edges of the network, a vertex can only be part of the refueling set if it is an endpoint. Therefore, $x$ must be in the interior of edge $(v_a, v_b)$, and has at least one endpoint on each side. Let $w_1$ and $w_2$ be the two closest endpoints to $x$ on each side. Since there are no vertices between $x$ and any of the two endpoints, all flows covered by $x$ are also covered by both endpoints.

Therefore, $T(x) \subseteq T(w_1) \cap T(w_2)$, which means that $T(x) \subseteq T(w_1)$ and $T(x) \subseteq T(w_2)$. Consequently, let $X^{*\prime}$ be an optimal solution that contains non-endpoints, $X^{*\prime} \nsubseteq EP$. Each non-endpoint, $x \in X^{*\prime}$, can be replaced by an endpoint, $w \in EP$, that covers at least the same set of routes, $T(x) \subseteq T(w)$, to get an optimal solution $X^* \subseteq EP$. □

Next, Theorem 2 derives an upper bound for the number of endpoints proving that the set of endpoints is finite.

**Theorem 2.** The ES algorithm for dedicated routes creates at most $4he$ endpoints, where $h = |H|$, $e = |E|$.

**Proof.** The ES algorithm attempts to establish a refueling set for each route $U_t \in H$ on each edge $(v_a, v_b) \in E$. If edge $(v_a, v_b) \in U_t$, then the entire $C(v_a, v_b)$ belongs to the refueling set, and the two vertices $v_a$ and $v_b$ are the endpoints. If edge $(v_a, v_b) \notin U_t$, then vehicles may deviate from their dedicated route to access this edge, subject to Constraints (6) and (7). A vehicle may enter the edge from either vertex $v_a$ or $v_b$, and for each entry vertex a refueling segment $RS^{(q)}(U_t; v_a, v_b)$, $q \in \{a, b\}$, is defined based on the remaining allowed deviation distance $\delta^{(q)}(U_t; v_a, v_b)$. If this distance is negative, $\delta^{(q)}(U_t; v_a, v_b) < 0$, vehicles dedicated to $U_t$ cannot reach vertex $v_q$ and no endpoints are added from this direction. If it equals zero, $\delta^{(q)}(U_t; v_a, v_b) = 0$, only the entry vertex qualifies as an endpoint. And if it is positive, $\delta^{(q)}(U_t; v_a, v_b) > 0$, the refueling segment begins at the entry vertex and extends a positive distance along the edge, resulting in two endpoints. Therefore, under the upper bound scenario, where the remaining allowed deviation distance is positive for both entry vertices, but the two refueling segments do not overlap, an edge may have up to four endpoints with respect to a single route. Since each of the $h$ routes can contribute up to four endpoints on each of the $e$ edges, the total number of endpoints across the entire network is bounded from above by $4he$. □

Based on Theorem 1, for any point located either at a vertex in $V$ or along an edge in $E$ there exists an endpoint that covers at least the same set of routes. Hence, for any solution that includes locations outside the set of endpoints we can always replace these locations with endpoints from $EP$ to obtain a solution that covers at least the same set of routes. This makes the set of endpoints a dominating set. Furthermore, Theorem 2 confirms that the set of endpoints is finite. Therefore, the set of endpoints constitutes a FDS.

Thus, since the set of endpoints forms a FDS, it is always possible to find an optimal solution to the dedicated routes refueling station location problem by establishing the set of endpoints, $EP$, then using these endpoints as candidate locations in a set covering model. The effort to establish the set of endpoints can be reduced by only considering edges $(v_a, v_b) \in E$ where at least one of the two entry vertices $v_a$, $v_b$ is within distance $D$ from at least one of the route's vertices. Take the network in Figure 8 as an example.

If we are looking for refueling sets that can cover route $U_1$, there is no need to consider any of the edges in $\{U_2 \cup (v_b, v_{12})\}$. This is because the distance from any potential deviation vertex in $V(U_1)$ to any potential entry vertex for these edges will exceed $D$; hence, violating constraint (7). Therefore, a set of potential refueling edges can be established for each route, $U_t \in H$. Let this set be $M(U_t)$ which is defined as follows:

$$M(U_t) = \{(v_a, v_b) | (v_a, v_b) \in U_t \; or \; \min\{l(v^d_{U_t}(v_a), v_a), l(v^d_{U_t}(v_b), v_b)\} \leq D\}. \qquad (26)$$

Using constraint (6) to reduce the size of $M(U_t)$ is more challenging. This is because satisfying this constraint depends on the specific network structure and the availability of deviation paths that keep the total travel distance below $R$. Figure 9 shows an example where the length of the route equals the driving range, $l(U_1) = R = 20$. A possible deviation path in this example can go from $v_3$ to $v_4$ and return to the route via $v_{10}$. This path will not increase the total trip length beyond $R$; therefore, a refueling station on edge $(v_3, v_4)$ or edge $(v_4, v_{10})$ can feasibly cover route $U_1$. However, if edge $(v_4, v_{10})$ does not exist, then edge $(v_3, v_4)$ cannot hold a refueling station that feasibly cover route $U_1$. This is because the total trip length will increase beyond $R$.

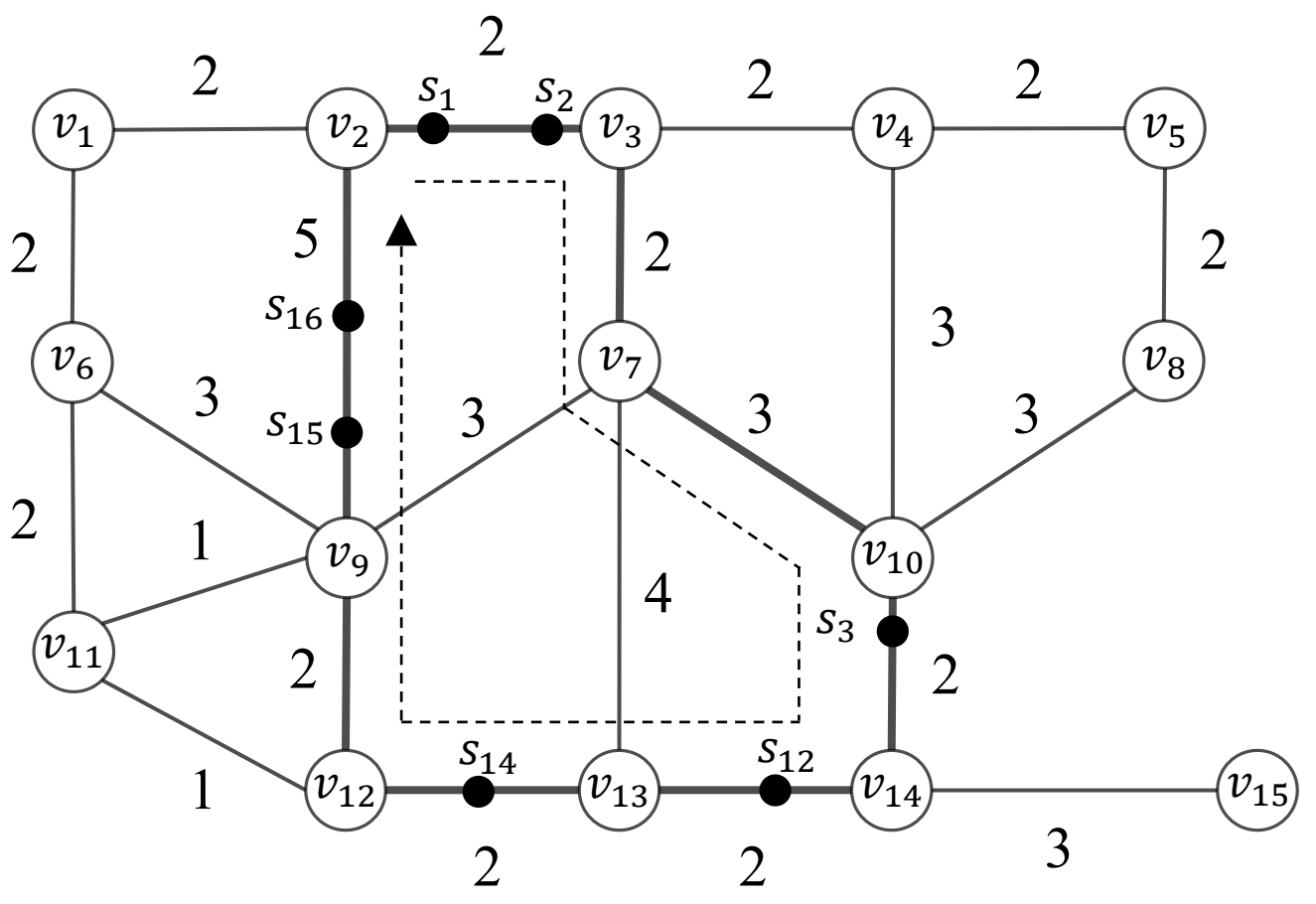

$U_1 = \{(v_2, v_3), (v_3, v_7), (v_7, v_{10}), (v_{10}, v_{14}), (v_{14}, v_{13}), (v_{13}, v_{12}), (v_{12}, v_9), (v_9, v_2)\}$, $l(U_1) = 20, D = 5, R = 20$.

Figure 9: Computational effort reduction

*4.3. Exact Algorithm for the Dedicated Routes Refueling Station Location Problem*

In this subsection, an exact polynomial-time algorithm is developed to find the refueling sets capable of covering each dedicated route in a given network and the overall set of endpoints. First, the algorithm establishes the set of routes to be covered, $H$. Then, for each route $U_t \in H$, we define the set of potential refueling edges, $M(U_t)$. After that, the algorithm establishes the refueling set capable of covering route $U_t$ on each edge in $M(U_t)$, and finds the set of endpoints. Finally, the endpoints are consolidated into the overall set of endpoints which will be used in a set covering model to find an optimal solution.

*ES Algorithm for Dedicated Routes*

Step 1: Generate the set of dedicated routes with positive flows that can be covered using a single refueling station, $H$.

$$H = \{U_t | l(U_t) \leq R, f(U_t) > 0\}.$$

Step 2: For each route $U_t \in H$, establish the set of potential refueling edges $M(U_t)$. Then find the refueling set capable of covering $U_t$ and its endpoints on each potential refueling edge in $M(U_t)$.

Sub-step 2.1: Establish the set $M(U_t)$ by including all edges in $U_t$. Then, for each vertex $v_i \in V(U_t)$ find the set of vertices within distance $D$ and add all its incident edges to $M(U_t)$. The set $M(U_t)$ becomes:

$$M(U_t) = \left\{(v_a, v_b) | (v_a, v_b) \in U_t \text{ or } \min\left\{l\left(v_{U_t}^d(v_a), v_a\right), l\left(v_{U_t}^d(v_b), v_b\right)\right\} \leq D\right\},$$

Where $v_{U_t}^d(v_a)$ and $v_{U_t}^d(v_b)$ are the deviation vertices in $V(U_t)$ to reach $v_a$ and $v_b$ respectively.

Sub-step 2.2: For each edge $(v_a, v_b) \in M(U_t)$, $(v_a, v_b) \notin U_t$, find $\beta^{(q)}(U_t; v_a, v_b)$ and $\delta^{(q)}(U_t; v_a, v_b)$ as follows:

$$\beta^{(q)}(U_t; v_a, v_b) =$$

$$\begin{cases} l(v_a, v_b), R - \begin{pmatrix} l(v_{U_t}^d(v_q), v_q) + l(v_w, v_{U_t}^b(s_{U_t}^f(s_{U_t}^l(v_q)))) + \\ l'(U_t(v_{U_t}^b(s_{U_t}^f(s_{U_t}^l(v_q))), v_{U_t}^d(v_q))) \end{pmatrix} \geq l(v_a, v_b), \\ \left( R - \begin{pmatrix} l(v_{U_t}^d(v_q), v_q) + l(v_q, v_{U_t}^b(s_{U_t}^f(s_{U_t}^l(v_q)))) + \\ l'(U_t(v_{U_t}^b(s_{U_t}^f(s_{U_t}^l(v_q))), v_{U_t}^d(v_q))) \end{pmatrix} \right) / 2, \text{Otherwise}, \end{cases}$$

$$q \in \{a, b\}; w \in \{a, b\}; w \neq q,$$

$$\delta^{(q)}(U_t; v_a, v_b) = \min\left\{D - l\left(v_{U_t}^d(v_q), v_q\right), \beta^{(q)}(U_t; v_a, v_b)\right\}, q \in \{a, b\}.$$

Sub-step 2.3: Determine the set of refueling segments that cover $f(U_t)$ on $(v_a, v_b) \in M(U_t)$ as follows:

If $(v_a, v_b) \in U_t$,

$$RS(U_t; v_a, v_b) = \{x | x \in C(v_a, v_b), (v_a, v_b) \in U_t\},$$

otherwise,

$$RS^{(q)}(U_t; v_a, v_b) = \{x | x \in C(v_a, v_b), (v_a, v_b) \notin U_t, l(v_q, x) \leq \delta^{(q)}(U_t; v_a, v_b)\}, q \in \{a, b\},$$

$$RS(U_t; v_a, v_b) = RS^{(a)}(U_t; v_a, v_b) \cup RS^{(b)}(U_t; v_a, v_b).$$

Sub-step 2.4: Determine the set of endpoints that cover $f(U_t)$ on $(v_a, v_b) \in M(U_t)$.

$$EP(U_t; v_a, v_b) = \{w^k_{U_t;v_a,v_b} | w^k_{U_t;v_a,v_b}, k = 1, \ldots, 4, \text{are endpoints of } RS(U_t; v_a, v_b)\}.$$

Step 3: Combine the sets of endpoints on each edge $(v_a, v_b) \in E$. Then establish the overall set of endpoints in the network as follows:

$$EP(v_a, v_b) = \cup_{U_t \in H} EP(U_t; v_a, v_b),$$

$$EP = \cup_{(v_a,v_b) \in E} EP(v_a, v_b).$$

Step 4: Determine the set of routes covered by each endpoint, $w \in EP$, and its total covered flow using:

$$T(w) = \{U_t | U_t \in H, w \in R(U_t; v_a, v_b), (v_a, v_b) \in E\},$$

$$F(w) = \sum_{U_t \in T(w)} f(U_t).$$

Next, we introduce Theorem 2 to quantify the required computational effort to solve the refueling station location problem for dedicated routes using our proposed algorithm.

**Theorem 3.** The computational complexity of the ES algorithm for dedicated routes is $O(he(n^2 + h))$, where $h = |H|$, $e = |E|$, and $n = |V|$.

**Proof.** Distances between all pairs of vertices in $V$ can be found using Johnson's algorithm in $O(n^2 \, log \, n + ne)$ [50]. Step 1 starts with calculating the lengths of all routes in the network. For any give route, the maximum number of unique edges can be the total number of edges in the network $e$. Additionally, each edge could appear multiple times in the route. Therefore, to calculate the length of the route, first the length of each edge will be multiplied by the number of times the edge appears in the route resulting in at most $e$ multiplications. Then these lengths are added in at most $e - 1$ additions. Hence at most $2e - 1$ operations are needed to calculate the length of one route. For all routes, at most $h2e - 1$ operations are needed. After that, two comparison operations are performed for each route, one to ensure that the length of the route is less than or equal to $R$, the other to ensure a positive flow. For all routes, $2h$ comparisons are performed. Therefore, Step 1 takes at most $h2e + 2h - 1$ operations. In Step 2, for each route, we establish the set of refueling edges and then find the set of endpoints on each edge that can cover the route. Sub-step 2.1 compares the distances between each vertex in the route and all other vertices in the network with the maximum allowed deviation distance, $D$. This takes at most $n$ comparisons for each vertex in the route, leading to an upper bound of $n^2/4$ comparisons for all vertices. After that, we include all incident edges of the vertices within $D$ from the route. Therefore, for each route there are $en^2$ union operations and a total of at most $hen^2$ operations for all routes. Sub-step 2.2 finds the remaining possible deviation distance after reaching an edge, $(v_a, v_b) \in M(U_t)$, $(v_a, v_b) \notin U_t$. In this sub-step, we perform at most 10 addition, multiplication, and comparison operations for each edge for a total of $10e$ operations. In Sub-steps 2.3 and 2.4, for each edge a series of comparisons will take place to determine the refueling segment and its endpoints for at most $5e$ and $4e$ operations, respectively. Considering the longest sub-step in Step 2, this step takes $hen^2$. In Step 3, first we combine the sets of endpoints on each edge. For each edge, there can be at most one set of endpoints per route for a total of $h$ sets. Therefore, for all edges there are $he$ union operations. After that, we combine the sets of endpoints on all edges into the set $EP$ which takes $e$ union operations. Hence, Step 3 takes $he + e$ union operations. Finally, in Step 4 we find the set of routes covered by each endpoint and its total covered flow. Based on Theorem 2, there can be at most $4he$ endpoints. For each endpoint on a given edge, all routes need to be considered to check if the endpoint belongs to a refueling segment of the route on that edge, this takes $4h^2e$. Then the flows of the routes that can be covered by each endpoint are added which takes $4h^2e$. Therefore, Step 4 takes $8h^2e$. Based on that and ignoring constant values, the dominant steps are Step 2 taking $hen^2$ and Step 4 taking $h^2e$. Therefore, the overall algorithm complexity is $O(he(n^2 + h))$. □

Based on Theorem 3, the computational complexity depends on the number of vertices, edges, and routes in the network. Theoretically, the number of routes can be extremely large given that every combination of 2 or more vertices and their connecting edges can be a potential route. However, in reality, the number of routes in any given city will be far less than the theoretical upper bound. We will discuss this further with practical examples in Section 5.

### *4.4. Mathematical Model*

Based on Theorem 1, there is always an optimal solution to the dedicated routes refueling station location problem where all stations are located at endpoints, $EP$. Therefore, after implementing the ES algorithm, the set of endpoints can be used in a set covering model to find an optimal solution. However, to reduce the size of the set covering problem, from every group of endpoints that covers the same set of routes, say $W = \{w_1, w_2 \ldots\}$, where $T(w_1) = T(w_2) = \cdots$, we can include only one endpoint, say $w_k \in W$, in the set of candidate locations. This will not affect the optimal solution since the set covering model does not look at the point's exact location, but the set of routes or the amount of flow covered by the candidate location. Therefore, let us define $C$ as a subset of endpoints, $C \subseteq EP$, that cover unique sets of routes. Let $X^*$ be an optimal solution to the set covering model and $T(X^*)$ is the set of routes covered by this solution, $T(X^*) = \bigcup_{w \in X^*} T(w)$. Now, if endpoint $w_k \in X^*$, we can replace $w_k$ by any point $x$ as long as $\{T(X^*) \backslash T(w_k)\} \cup T(x) = T(X^*)$ to find an alternative optimal solution. Below, we present a set covering model to find the minimum number of refueling stations required to cover all routes in the network.

$C$: Set of candidate locations, $|C| = c$.

$H$: Set of routes, $|H| = h$.

$A$: $(c \times h)$ binary matrix where element $a_{w,U_t} = 1$ if route $U_t \in H$ can be covered by a refueling station in location $w \in C$.

$x_w$: Decision variable, $x_w = \begin{cases} 1, \text{if a refueling station is located at } w \in C, \\ 0, \text{otherwise.} \end{cases}$

Set covering model:

$$\text{Minimize: } p = \sum_{w \in C} x_w, \tag{27}$$

Subject to:

$$\sum_{w \in C} a_{w,U_t} x_w \geq 1, \forall U_t \in H, \tag{28}$$

$$x_w \in \{0,1\}, w \in C. \tag{29}$$

As mentioned above, after solving the set covering model, any selected endpoint in the optimal solution, $w_k \in C$, can be replaced by a point $x$, as long as $x$ covers all routes that are uniquely covered by $w_k$, i.e., $\{T(X^*) \backslash T(w_k)\} \cup T(x) = T(X^*)$ to find an alternative optimal solution. Lemma 2 helps us find interior points that can be included in the optimal solution.

**Lemma 2:** Let $w_a$, $w_b$ be two adjacent endpoints in $EP$ that belong to an edge $(v_a, v_b) \in E$, $w_a \neq w_b$. All interior points of the line segment connecting $w_a$ and $w_b$, int $(C(w_a, w_b))$, cover the same set of routes. Moreover, let $x \in \text{int}\big(C(w_a, w_b)\big)$, $T(x) \subseteq T(w_a) \cup T(w_b)$.

**Proof:** Without loss of generality, let $l(v_a, w_a) < l(v_a, w_b)$ and $l(v_b, w_b) < l(v_b, w_a)$. Dedicated route $U_t \in T(w_q)$, $q \in \{a, b\}$, must have a refueling segment on edge $(v_a, v_b)$ that starts at $v_q$ and includes $w_q$, $RS^{(q)}(U_t; v_a, v_b)$, $q \in \{a, b\}$. $RS^{(q)}(U_t; v_a, v_b)$ could end at $w_q$. In this case, vehicles dedicated to route $U_t$ and entering edge $(v_a, v_b)$ using vertex $v_q$, $q \in \{a, b\}$, cannot reach any point in int $(C(w_a, w_b))$. Or $RS^{(q)}(U_t; v_a, v_b)$ could end at $w_r$, $q \in \{a, b\}$, $r \in \{a, b\}$, $r \neq q$. In this case, vehicles dedicated to route $U_t$ and entering edge $(v_a, v_b)$ using vertex $v_q$, can reach all points in $C(w_a, w_b)$. $RS^{(q)}(U_t; v_a, v_b)$ cannot end at any point in int $\big(C(w_a, w_b)\big)$ since $w_a$ and $w_b$ are adjacent endpoints. Additionally, given that a vertex can only be part of a refueling segment if it is an endpoint, int $(C(w_a, w_b))$ cannot have any vertices.

Therefore, vehicles reaching any interior point in $\text{int}\,(C(w_a, w_b))$ must have come through one of the endpoints and cannot exit edge $(v_a, v_b)$ before reaching the other endpoint. Based on that, all interior points of line segment $C(w_a, w_b)$ cover the same set of routes and $T(x) \subseteq T(w_a) \cup T(w_b)$, where $x \in \text{int}\left(C(w_a, w_b)\right)$. □

Based on Lemma 2, the problem can be solved optimally considering $EP$ as the set of candidate locations. Then, for any two adjacent endpoints $w_a, w_b \in EP$ it is sufficient to examine a single point $x$, $x \in \text{int}\left(C(w_a, w_b)\right)$, to know the set of routes that can be covered by any point in $\text{int}\left(C(w_a, w_b)\right)$. Then we can determine if any of these points can replace an endpoint in the current optimal solution to create an alternative optimal solution. Finally, it is important to note that $T(x) \neq T(w_a) \cap T(w_b)$. Figure 10 shows an example where vehicles of route $U_1$ can reach the two adjacent endpoints $w_a$ and $w_b$ from the two sides of edge $(v_a, v_b)$, but cannot reach any interior point in $\text{int}\left(C(w_a, w_b)\right)$. Therefore, $T(x) \neq T(w_a) \cap T(w_b)$, however $T(x) \subseteq T(w_a) \cup T(w_b)$.

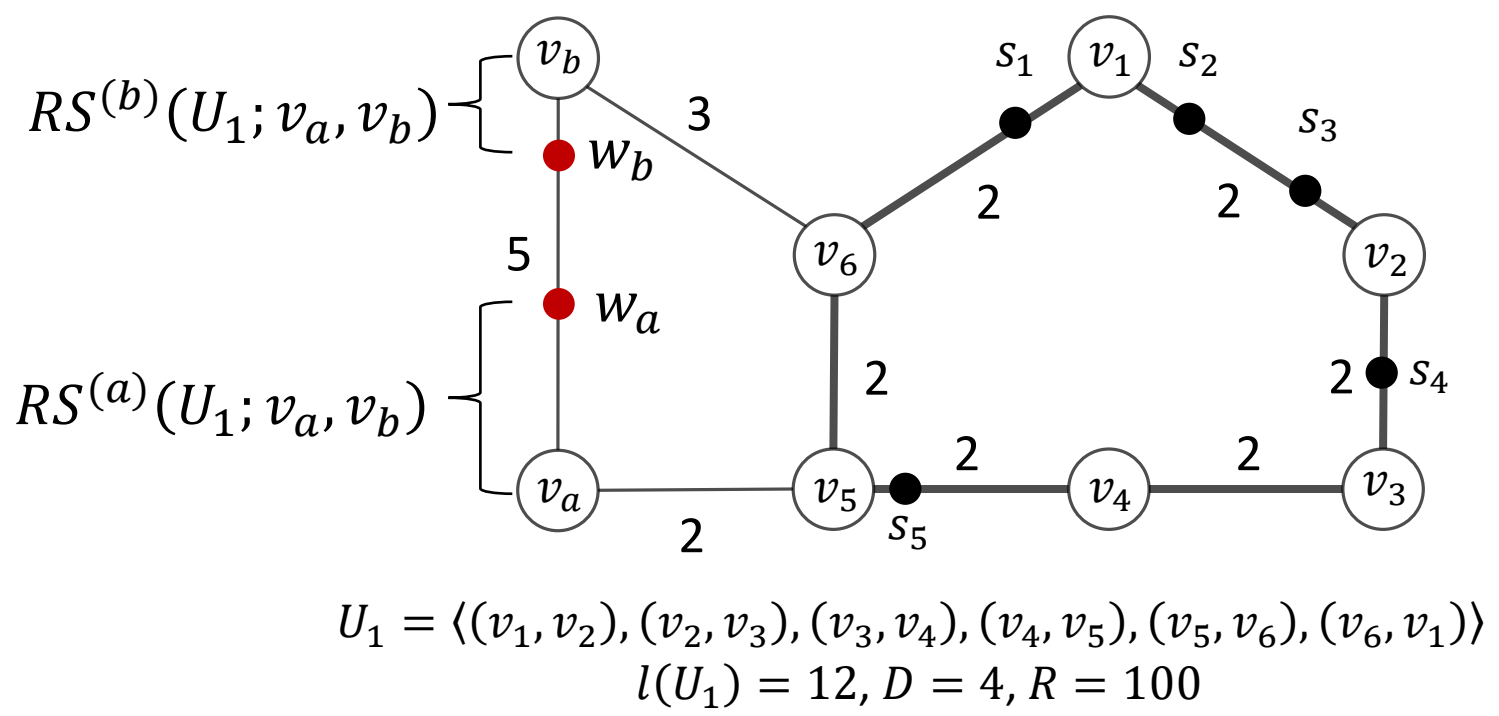

Figure 10: Interior points coverage example

## 5. Numerical Experiments

To illustrate the proposed methodology, we applied it to two well-known and realistic transportation networks: the Sioux Falls network and the San Antonio network. The Sioux Falls network [51], shown in Figure 11, is widely used as a benchmark in transportation research. We use this network to illustrate the detailed steps of our methodology. The San Antonio network is larger and more complex, representing a practical urban scenario. It is used to provide insights into the scalability of the proposed methodology. This section is divided into four subsections: Sioux Falls Network, San Antonio Network, Sensitivity Analysis, and Computational Effort.

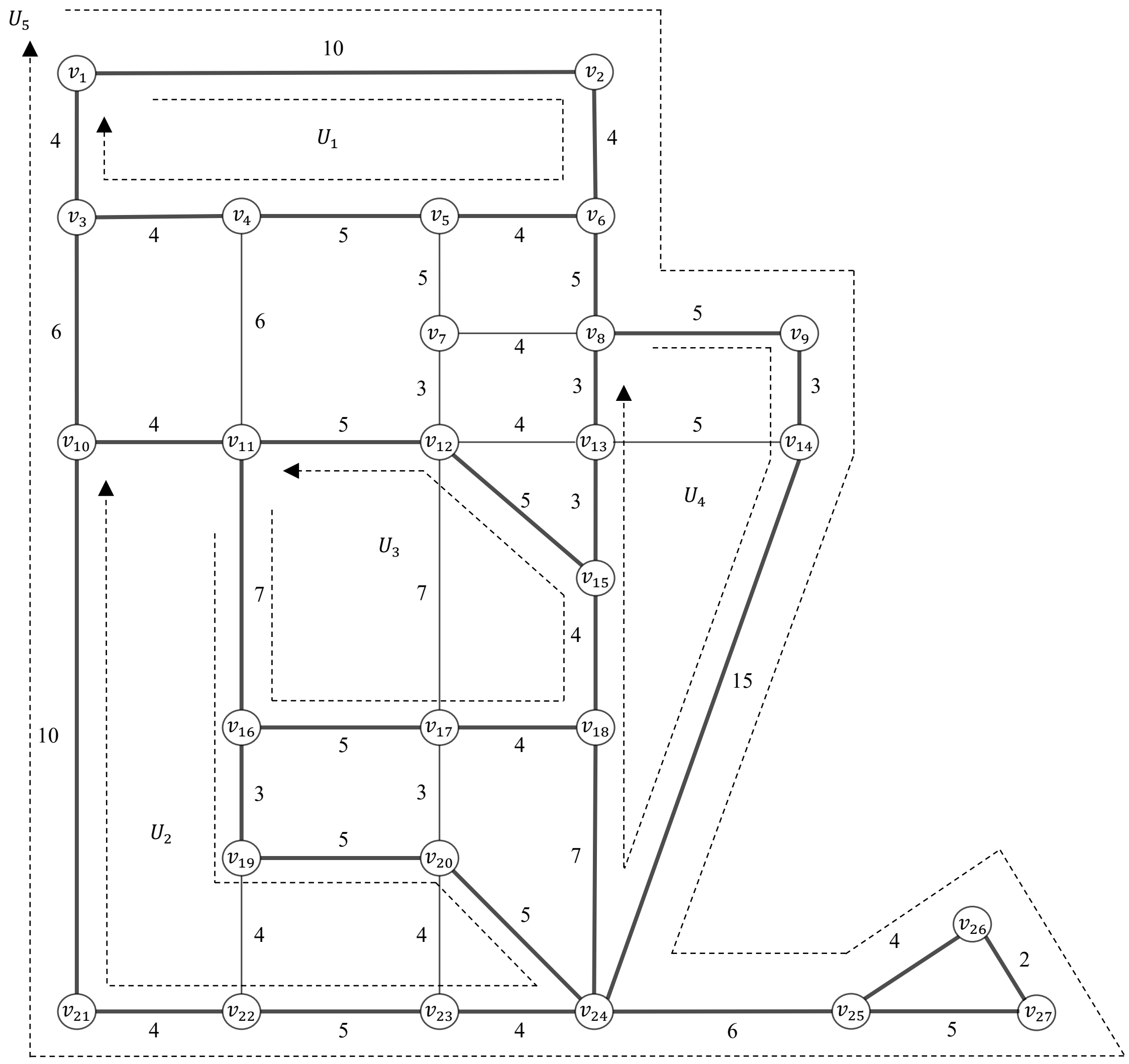

Figure 11: Dedicated routes example for Sioux Falls network

### *5.1. Sioux Falls Network*

The Sioux Falls network consists of 27 nodes, 40 edges, and 5 dedicated routes. The dedicated routes are listed in Table 1 and they are presented in Figure 11 using dashed arrows. Our network has undirected roads and a stop just before and just after every vertex. We assume that all vehicles have a fixed driving range of $R = 100$ miles, and a maximum allowed deviation distance toward a refueling station $D = 4$ miles.

For illustration, let us consider edge $(v_6, v_8)$ for finding refueling sets. This edge is part of route $U_5$, therefore any point on the edge covers route $U_5$. The shortest distance from each of vertices $v_6$ and $v_8$ to routes $U_2$ and $U_3$ is greater than $D = 4$. Therefore, these routes cannot be covered on edge $(v_6, v_8)$. Next, consider route $U_4$:

$$\beta^{(6)}(U_4; v_6, v_8) = \begin{cases} l(v_6, v_8), R - \Big(l(v_8, v_6) + l(v_8, v_8) + l'\big(U_4(v_8, v_8)\big)\Big) \geq l(v_6, v_8), \\ \Big(R - \Big(l(v_8, v_6) + l(v_6, v_8) + l'\big(U_4(v_8, v_8)\big)\Big)\Big)/2, \text{Otherwise.} \end{cases}$$

$$= \begin{cases} 5, 100 - (5 + 0 + 40) \geq 5, \\ \big(100 - (5 + 5 + 40)\big)/2, \text{Otherwise.} \end{cases}$$

Therefore $\beta^{(6)}(U_4; v_6, v_8) = 5$. Similarly, $\beta^{(8)}(U_4; v_6, v_8) = 5$. Now, the remaining deviation distances for a vehicle dedicated to route $U_4$ when it reaches vertices $v_6$ and $v_8$ are calculated as follows:

$\delta^{(6)}(U_4; v_6, v_8) = \min\{D - l(v_8, v_6), \beta^{(6)}(U_4; v_6, v_8)\} = \min\{-1, 5\} = -1$,

$\delta^{(8)}(U_4; v_6, v_8) = \min\{D - l(v_8, v_8), \beta^{(8)}(U_4; v_6, v_8)\} = \min\{4, 5\} = 4$.

Table 1: Dedicated routes for the Sioux Falls network

| Route | Sequence of edges | Length | Flow |
|---|---|---|---|
| $U_1$ | $\langle(v_1, v_2), (v_2, v_6), (v_6, v_5), (v_5, v_4), (v_4, v_3), (v_3, v_1)\rangle$ | 31 | 25 |
| $U_2$ | $\langle(v_{10}, v_{11}), (v_{11}, v_{16}), (v_{16}, v_{19}), (v_{19}, v_{20}), (v_{20}, v_{24}), (v_{24}, v_{23}),$ $(v_{23}, v_{22}), (v_{22}, v_{21}), (v_{21}, v_{10})\rangle$ | 47 | 30 |
| $U_3$ | $\langle(v_{11}, v_{16}), (v_{16}, v_{17}), (v_{17}, v_{18}), (v_{18}, v_{15}), (v_{15}, v_{12})\ (v_{12}, v_{11})\rangle$ | 30 | 15 |
| $U_4$ | $\langle(v_8, v_9), (v_9, v_{14}), (v_{14}, v_{24})\ (v_{24}, v_{18}), (v_{18}, v_{15}), (v_{15}, v_{13}),\ (v_{13}, v_8)\rangle$ | 40 | 10 |
| $U_5$ | $\langle(v_1, v_2), (v_2, v_6), (v_6, v_8), (v_8, v_9), (v_9, v_{14}), (v_{14}, v_{24}), (v_{24}, v_{25}),$ $(v_{25}, v_{26}), (v_{26}, v_{27}), (v_{27}, v_{25}), (v_{25}, v_{24}), (v_{24}, v_{23}), (v_{23}, v_{22}),$ $(v_{22}, v_{21}), (v_{21}, v_{10}), (v_{10}, v_3), (v_3, v_1)\rangle$ | 98 | 5 |

Since $\delta^{(6)}(U_4; v_6, v_8)$ is negative, vehicles dedicated to route $U_4$ cannot reach edge $(v_6, v_8)$ via vertex $v_6$, $RS^{(6)}(U_4; v_6, v_8) = \phi$. However, they can reach this edge via vertex $v_8$ and have 4 miles remaining of the maximum allowed deviation distance, $RS^{(6)}(U_4; v_5, v_6) = \{x | x \in C(v_5, v_6), l(v_6, x) \leq 4\}$. Next, consider route $U_1$. The values for $\beta^{(q)}(U_1; v_6, v_8)$, $q \in \{6,8\}$, are calculated as follows:

$$\beta^{(6)}(U_1; v_6, v_8) = \begin{cases} l(v_6, v_8), R - \Big(d(v_6, v_6) + d(v_8, v_6) + l'\big(U_1(v_6, v_6)\big)\Big) \geq l(v_6, v_8), \\ \Big(R - \Big(d(v_6, v_6) + d(v_6, v_8) + l'\big(U_1(v_6, v_6)\big)\Big)\Big)/2, \ \text{Otherwise.} \end{cases}$$

$$= \begin{cases} 5, 100 - (0 + 5 + 31) \geq 5, \\ \big(100 - (0 + 5 + 31)\big)/2, \text{Otherwise.} \end{cases}$$

Therefore $\beta^{(6)}(U_1; v_6, v_8) = 5$. Similarly, $\beta^{(8)}(U_1; v_6, v_8) = 5$. Now, the remaining deviation distances for a vehicle dedicated to route $U_1$ when it reaches vertices $v_6$ and $v_8$ are calculated as follows:

$\delta^{(6)}(U_1; v_6, v_8) = \min\{4 - d(v_6, v_6), \beta^{(6)}(U_1; v_6, v_8)\} = \min\{4, 5\} = 4$,

$\delta^{(8)}(U_1; v_6, v_8) = \min\{4 - d(v_6, v_8), \beta^{(6)}(U_1; v_6, v_8)\} = \min\{-1, 5\} = -1$.

Since $\delta^{(8)}(U_1; v_6, v_8)$ is negative, then vehicles dedicated to route $U_1$ cannot reach edge $(v_6, v_8)$ via vertex $v_8$, $RS^{(8)}(U_1; v_6, v_8) = \phi$. However, they can reach this edge via vertex $v_6$ and have 4 miles remaining of the maximum allowed deviation distance, $RS^{(6)}(U_1; v_6, v_8) = \{x | x \in C(v_6, v_8), l(v_6, x) \leq 4\}$. The refueling sets and endpoints on edge $(v_6, v_8)$ are shown in Figure 12.

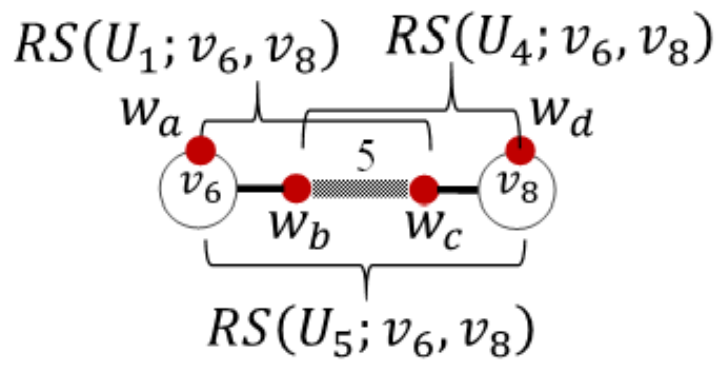

Figure 12: Refueling sets on edge $(v_6, v_8)$

Note that any point in the line segment $C(w_b, w_c)$ can cover routes $U_1$, $U_4$, and $U_5$. Similarly, we implemented the ES algorithm for dedicated routes to all edges in the network. The complete set of endpoints is listed in Table 2 and shown in Figure 13.

Table 2: List of endpoints and their covered routes for the Sioux Falls network

| Endpoints $(w)$ | Set of covered routes $T(w)$ | Total flow covered $F(w)$ in round trips/day |
|---|---|---|
| $w_9, w_{10}$ | $U_1$ | 25 |
| $w_{46}$ | $U_3$ | 15 |
| $w_{65}, w_{66}, w_{67}$ | $U_5$ | 5 |
| $w_2, w_{31}$ | $U_1, U_2, U_3$ | 70 |
| $w_1, w_{32}$ | $U_1, U_2, U_5$ | 60 |
| $w_3$ | $U_1, U_3$ | 40 |
| $w_4, w_{11}$ | $U_1, U_4, U_5$ | 40 |
| $w_5, w_6, w_7, w_8, w_{57}, w_{64}$ | $U_1, U_5$ | 30 |
| $w_{30}, w_{33}, w_{34}, w_{36}, w_{40}, w_{41}, w_{42}, w_{43}, w_{50}, w_{52}, w_{53}, w_{54}, w_{55}, w_{56}$ | $U_2, U_3$ | 45 |
| $w_{22}, w_{26}, w_{35}, w_{37}, w_{51}$ | $U_2, U_3, U_4$ | 55 |
| $w_{44}, w_{63}$ | $U_2, U_3, U_5$ | 50 |
| $w_{24}, w_{25}, w_{27}, w_{38}, w_{59}, w_{60}$ | $U_2, U_4, U_5$ | 45 |
| $w_{39}, w_{45}, w_{61}, w_{62}$ | $U_2, U_5$ | 35 |
| $w_{12}, w_{16}, w_{17}, w_{18}, w_{20}, w_{21}, w_{23}, w_{29}, w_{47}, w_{48}, w_{49}$ | $U_3, U_4$ | 25 |
| $w_{13}, w_{14}, w_{15}, w_{19}, w_{28}, w_{58}$ | $U_4, U_5$ | 15 |

Now, we use the set covering model (27) - (29) to find the minimum number of refueling stations needed to cover the five routes and the optimal locations. Only one endpoint from each row in Table 2 needs to be included in the set of candidate locations $C$. Therefore, rather than having 67 candidate locations, only 15 are needed in the set covering model. We used the first endpoint from each row in Table 2 as a candidate location and solved the model using IBM ILOG CPLEX 12.10.0.

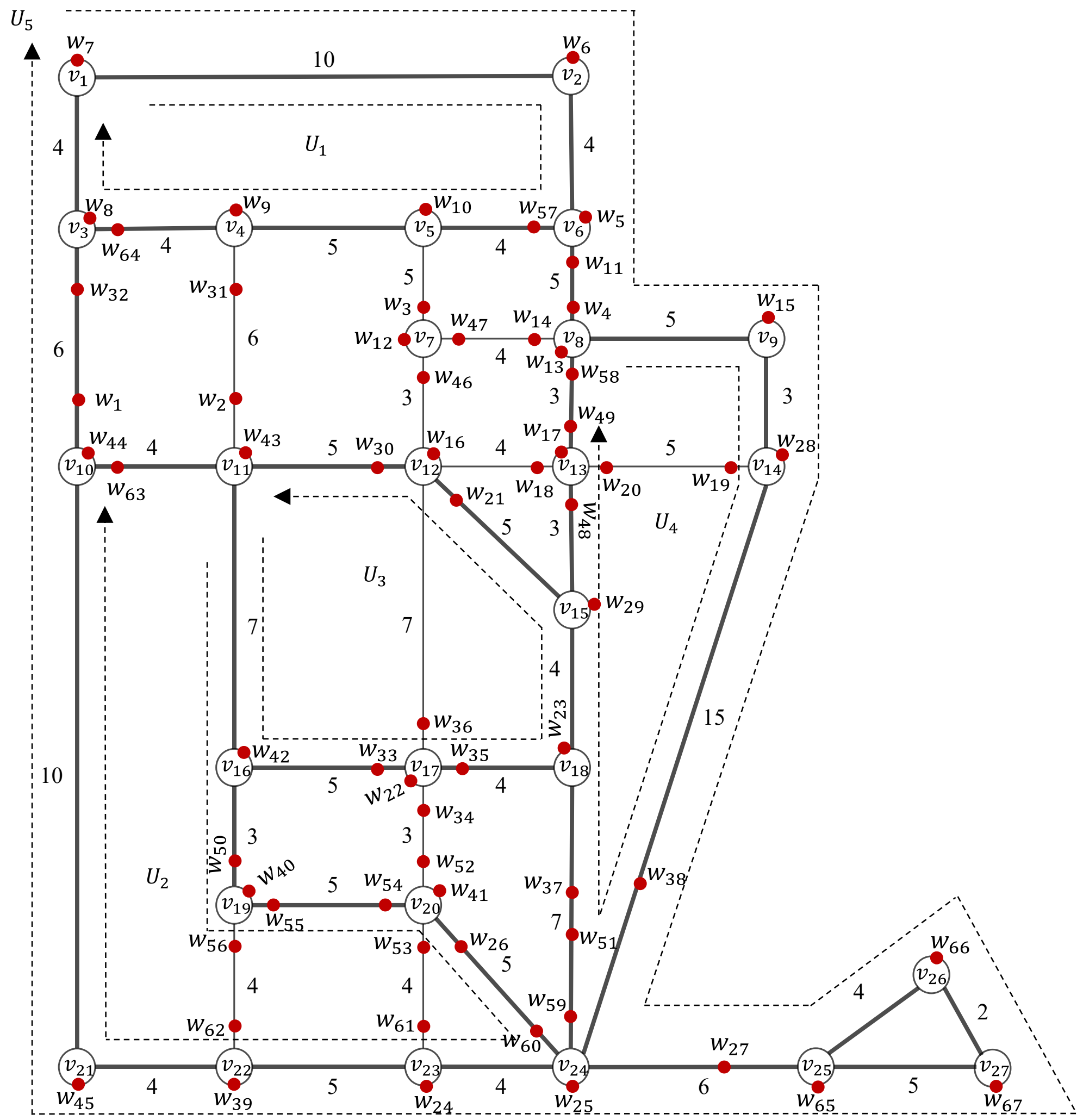

Figure 13: Dedicated routes example complete set of endpoints

The minimum number of refueling stations required to cover all routes in this example is two. One optimal solution is to locate a refueling station at $w_2$ to cover routes $U_1$, $U_2$, $U_3$, and a refueling station at $w_4$ to cover routes $U_1$, $U_4$, $U_5$. Note that $w_{31}$ covers the same set of routes as $w_2$, therefore we can create an alternative optimal solution by replacing $w_2$ in the optimal solution by $w_{31}$. Moreover, we can pick any point in int $C(w_2, w_{31})$, say the midpoint of edge $(v_4, v_{11})$, and determine the set of routes that can be covered by this point using our two coverage constraints (6) and (7). It is easy to verify that the range of $R = 100$ miles is enough to reach this point after completing at least one trip around any of the routes in this network except $U_5$. Also, this point is only 3 miles away from potential deviation vertices of routes $U_1$, $U_2$, and $U_3$ while it is 12 miles away from the potential deviation vertex of route $U_4$. Therefore, this point satisfies the two coverage constraints for routes $U_1$, $U_2$, and $U_3$. By Lemma 2, all points in int $C(w_2, w_{31})$ cover routes $U_1$, $U_2$, and $U_3$ and can replace $w_2$ in the optimal solution. In fact, any point that covers routes $U_2$ and $U_3$ can replace $w_2$ to find an alternative solution.

This flexibility extends to $w_4$ as well. Instead of placing a station at $w_4$ to cover routes $U_1$, $U_4$, and $U_5$, we can select $w_{11}$ as an alternative location. Moreover, any point between the endpoints of $w_4$ and $w_{11}$ can serve as an optimal location, allowing decision-makers to incorporate other selection criteria, such as land prices, accessibility, and expected population growth, when finalizing station placement.

In contrast, using the discrete approach with the set of vertices as candidate locations, the best solution requires three refueling stations to cover all routes. Table 3 shows the set of routes covered by each vertex in the network. One possible optimal solution under the discrete approach involves locating one refueling station at $v_{10}$, to cover routes $U_2$, $U_3$, and $U_5$, one refueling station at $v_8$ to cover $U_4$, and a third refueling station at $v_6$ to cover $U_1$. Compared to two refueling stations using our approach.

Table 3: List of vertices and their covered routes in the Sioux Falls network

| Vertices ($v_i$) | Set of covered routes $T(v_i)$ | Total flow covered $F(v_i)$ in round trips/day |
|---|---|---|
| $v_4$, $v_5$ | $U_1$ | 25 |
| $v_1$, $v_2$, $v_3$, $v_6$ | $U_1$, $U_5$ | 30 |
| $v_{19}$, $v_{20}$, $v_{16}$, $v_{11}$ | $U_2$, $U_3$ | 45 |
| $v_{17}$ | $U_2$, $U_3$, $U_4$ | 55 |
| $v_{10}$ | $U_2$, $U_3$, $U_5$ | 50 |
| $v_{23}$, $v_{24}$ | $U_2$, $U_4$, $U_5$ | 45 |
| $v_{21}$, $v_{22}$ | $U_2$, $U_5$ | 35 |
| $v_7$, $v_{12}$, $v_{13}$, $v_{18}$, $v_{15}$ | $U_3$, $U_4$ | 25 |
| $v_8$, $v_9$, $v_{14}$ | $U_4$, $U_5$ | 15 |
| $v_{25}$, $v_{26}$, $v_{27}$ | $U_5$ | 5 |

*5.2. San Antonio Network*

We evaluated the scalability and effectiveness of the ES algorithm by applying it to the San Antonio, Texas, network, a larger and more complex transportation system compared to the Sioux Falls network. Figure 14 shows a map of the major roads in San Antonio with 32 vertices and 52 edges. We selected a sample of six routes from the website of VIA Metropolitan Transit, a major public transportation company in San Antonio. Table 4 shows the lengths and flows of the selected routes.

Table 4: Dedicated routes for the San Antonio network

| Route | Sequence of edges | Length | Flow |
|---|---|---|---|
| $U_1$ | $\langle(v_5, v_3), (v_3, v_2), (v_2, v_{32}), (v_{32}, v_1), (v_1, v_{20}), (v_{20}, v_4), (v_4, v_5)\rangle$ | 52 | 15 |
| $U_2$ | $\langle(v_{18}, v_{19}), (v_{19}, v_4), (v_4, v_{20}), (v_{20}, v_{21}), (v_{21}, v_{18})\rangle$ | 43 | 20 |
| $U_3$ | $\langle(v_{16}, v_{17}), (v_{17}, v_{18}), (v_{18}, v_{22}), (v_{22}, v_{15}), (v_{15}, v_{12}), (v_{12}, v_{16})\rangle$ | 37 | 15 |
| $U_4$ | $\langle(v_{13}, v_{12}), (v_{12}, v_{15}), (v_{15}, v_{22})\ (v_{22}, v_{23}), (v_{23}, v_{14}), (v_{14}, v_{13})\rangle$ | 39 | 30 |
| $U_5$ | $\langle(v_{13}, v_{10}), (v_{10}, v_9), (v_9, v_8), (v_8, v_{11}), (v_{11}, v_{12}), (v_{12}, v_{13})\rangle$ | 44 | 25 |
| $U_6$ | $\langle(v_9, v_2), (v_2, v_3), (v_3, v_5), (v_5, v_6), (v_6, v_{11}), (v_{11}, v_8), (v_8, v_9)\rangle$ | 43 | 20 |

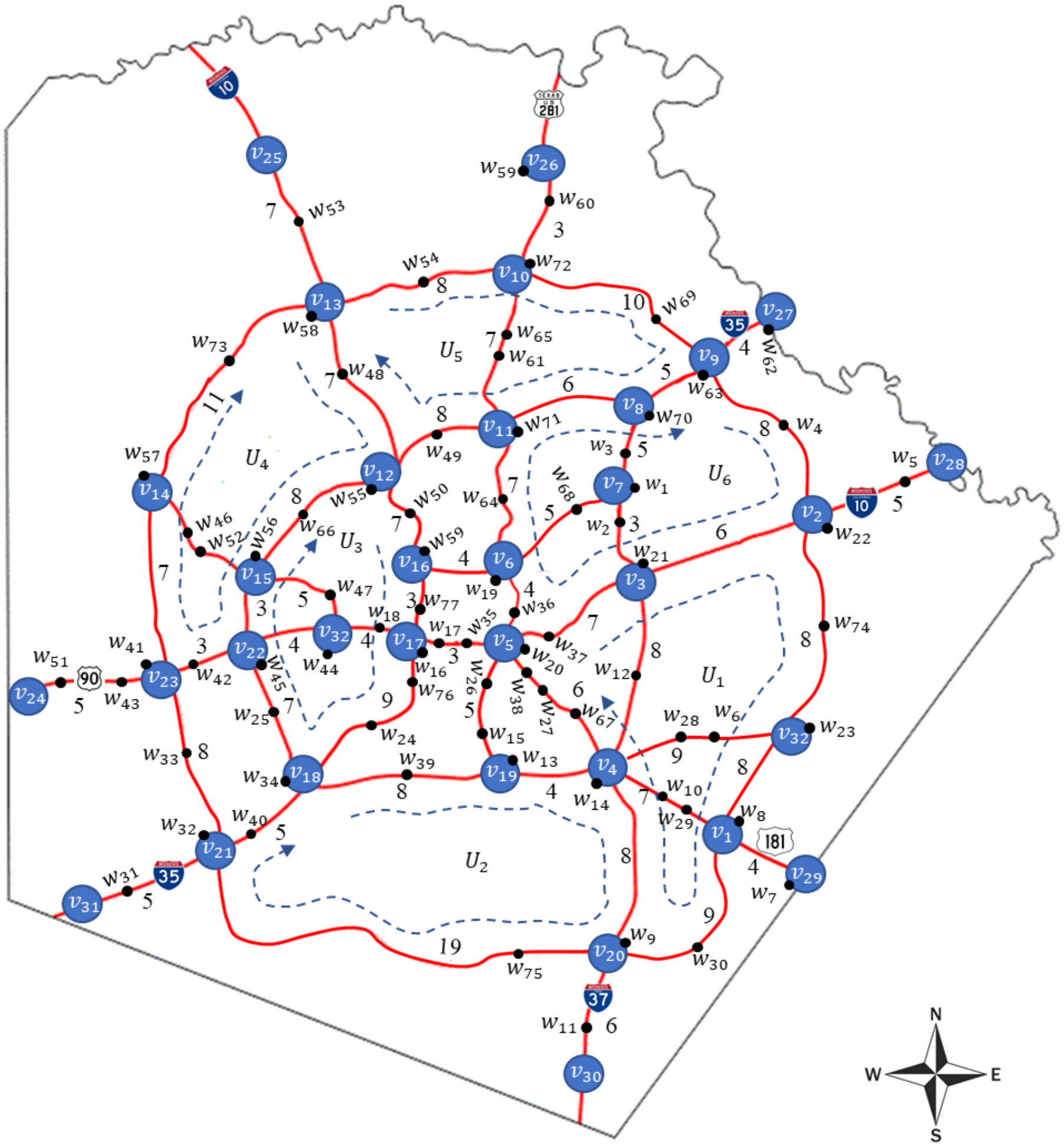

Figure 14: Dedicated routes example for San Antonio network

We assume that all vehicles have a fixed driving range of $R = 100$ miles, and a maximum allowed deviation distance toward a refueling station $D = 4$ miles. Applying the ES algorithm to the San Antonio network results in 77 endpoints, listed in Table 5. The analysis reveals that the minimum number of refueling stations required to cover all six routes is two. An optimal solution is to place a refueling station at $w_{26}$, which provides coverage for $U_1$, $U_2$, $U_3$, and $U_6$, and a refueling station at $w_{53}$, which covers $U_4$ and $U_5$. Additionally, $w_{54}$ can replace $w_{53}$ to create an alternative solution, offering flexibility in station placement.

Table 5: List of endpoints and their covered routes for the San Antonio network

| Endpoints ($w$) | Set of covered routes $T(w)$ | Total flow covered $F(w)$ in round trips/day |
|---|---|---|
| $w_6$, $w_7$, $w_8$, $w_{23}$ | $U_1$ | 15 |
| $w_{31}$, $w_{32}$ | $U_2$ | 20 |
| $w_{51}$ | $U_4$ | 30 |
| $w_{60}$, $w_{62}$, $w_{72}$ | $U_5$ | 25 |
| $w_9$, $w_{10}$, $w_{11}$, $w_{13}$, $w_{14}$, $w_{28}$, $w_{29}$, $w_{30}$, $w_{75}$ | $U_1$, $U_2$ | 35 |
| $w_{26}$ | $U_1$, $U_2$, $U_3$, $U_6$ | 70 |
| $w_{12}$, $w_{15}$, $w_{27}$, $w_{67}$ | $U_1$, $U_2$, $U_6$ | 55 |
| $w_{76}$, $w_{77}$ | $U_1$, $U_3$ | 30 |
| $w_{16}$, $w_{17}$, $w_{18}$, $w_{19}$, $w_{20}$, $w_{35}$, $w_{36}$, $w_{37}$, $w_{38}$ | $U_1$, $U_3$, $U_6$ | 50 |
| $w_3$, $w_4$ | $U_1, U_5, U_6$ | 60 |
| $w_1$, $w_2$, $w_5$, $w_{21}$, $w_{22}$, $w_{68}$, $w_{74}$ | $U_1$, $U_6$ | 35 |
| $w_{25}$ | $U_2$, $U_3$, $U_4$ | 65 |
| $w_{33}$ | $U_2$, $U_4$ | 50 |
| $w_{41}$, $w_{42}$, $w_{43}$, $w_{44}$, $w_{45}$, $w_{46}$, $w_{47}$, $w_{52}$, $w_{56}$, $w_{57}$ | $U_3$, $U_4$ | 45 |
| $w_{48}$, $w_{50}$, $w_{55}$, $w_{58}$, $w_{66}$, $w_{73}$ | $U_3$, $U_4$, $U_5$ | 70 |
| $w_{49}$ | $U_3$, $U_4$, $U_5$, $U_6$ | 90 |
| $w_{53}$, $w_{54}$ | $U_4$, $U_5$ | 55 |
| $w_{59}$, $w_{61}$, $w_{63}$, $w_{64}$, $w_{65}$, $w_{69}$, $w_{70}$, $w_{71}$ | $U_5$, $U_6$ | 45 |

Several alternative optimal solutions exist. For example, we can place a refueling station at $w_{49}$, which provides coverage for routes $U_3$, $U_4$, $U_5$, and $U_6$. Then to cover routes $U_1$ and $U_2$, a station can be located at any of the endpoints within the set $\{w_9, w_{10}, w_{11}, w_{13}, w_{14}, w_{28}, w_{29}, w_{30}, w_{75}\}$. Moreover, by Lemma 2, any point between $w_{10}$ and $w_{29}$ on edge $(v_1, v_4)$ can cover routes $U_1$ and $U_2$, offering further flexibility in station placement.

Similar to the Sioux Falls example, the best solution following the discrete approach and using the set of vertices as the set of candidate locations requires at least three refueling stations to cover all flows. From Table 6, an optimal station placement under the discrete approach involves locating a station at $v_5$, $v_6$, or $v_{17}$ to cover $U_1$, $U_3$, and $U_6$, another station at $v_{12}$ or $v_{13}$ to cover $U_3$, $U_4$, and $U_5$, and a third station at $v_{21}$ to cover $U_2$.

Table 6: List of vertices and their covered routes in the San Antonio network

| Vertices ($v_i$) | Set of covered routes $T(v_i)$ | Total flow covered $F(v_i)$ in round trips/day |
|---|---|---|
| $v_1, v_{29}, v_{32}$ | $U_1$ | 15 |
| $v_{21}$ | $U_2$ | 20 |
| $v_{10}, v_{27}$ | $U_5$ | 25 |
| $v_4, v_{19}, v_{20}$ | $U_1, U_2$ | 35 |
| $v_5, v_6, v_{17}$ | $U_1, U_3, U_6$ | 50 |
| $v_2, v_3, v_7$ | $U_1, U_6$ | 35 |
| $v_{14}, v_{15}, v_{22}, v_{23}, v_{32}$ | $U_3, U_4$ | 45 |
| $v_{12}, v_{13}$ | $U_3, U_4, U_5$ | 70 |
| $v_8, v_9, v_{11}, v_{26}$ | $U_5, U_6$ | 45 |

### *5.3. Sensitivity Analysis*

Sensitivity analysis is conducted to evaluate the robustness of the proposed ES algorithm in identifying a FDS. The analysis is performed on the Sioux Falls network, which has five dedicated routes. The objective is to assess how varying the deviation distance $D$ and the driving range $R$ affects the number and placement of refueling stations. To analyze the sensitivity of the model, we selected three different cases:

- Case 1: $D = 3$ miles and $R = 100$ miles
- Case 2: $D = 11$ miles and $R = 100$ miles
- Case 3: $D = 4$, miles and $R = 106$ miles

For comparison, the Sioux Falls network example presented in Section 5.1 is referred to as the original case throughout the sensitivity analysis. The results of Cases 1, 2, and 3 are compared against this original case.

When the maximum allowed deviation distance is small, such as in Case 1, vehicles have limited flexibility to deviate from their routes, which restricts the number of accessible edges outside the route itself. Consequently, a small $D$ reduces the possibility for a given endpoint to cover multiple routes which could lead to an increase in the number of refueling stations needed. Additionally, changing the value of $D$ can affect the number of endpoints in the network. For example, in the original case, $D = 4$ miles, the total number of endpoints was 67. However, when $D$ is reduced to 3 miles in Case 1, Figure A.1 (see appendix), the total number of endpoints decreased to 59. A notable example of this effect is observed when $D = 4$, miles and $R = 100$ miles, an intersection area emerges along edge $(v_3, v_{10})$ where the remaining deviation distances $\delta^{(3)}(U_1; v_3, v_{10})$ and $\delta^{(10)}(U_2; v_3, v_{10})$ equal 4 miles with corresponding endpoints $w_1$ and $w_{32}$, respectively. Therefore, there is an intersection segment, $(w_1, w_{32})$, with the length $l(w_1, w_{32}) = 2$ miles, that can cover both routes $U_1$ and $U_2$. Now, if we reduce the deviation distance to $D = 3$ miles, the remaining deviation distances $\delta^{(3)}(U_1; v_3, v_{10})$ and $\delta^{(10)}(U_2; v_3, v_{10})$ decrease to 3 miles and the two endpoints will coincide at the midpoint of edge $(v_3, v_{10})$. Consequently, the intersection segment will be reduced to a single point, that is the midpoint, and it will be the only location on edge $(v_3, v_{10})$ that can simultaneously cover both routes. Further reduction of the maximum allowed deviation distance (e.g., $D =$ 2 miles) will reduce $\delta^{(3)}(U_1; v_3, v_{10})$ and $\delta^{(10)}(U_2; v_3, v_{10})$ to 2 miles, where the endpoints will separate again, but there will be no intersection segment and none of the points on edge $(v_3, v_{10})$ can cover both routes.

Table 7 contains the list of endpoints for Case 1 along with the combination of routes that can be covered by each endpoint and the corresponding flows. Using this set of endpoints in the mathematical model yields an optimal solution using two refueling stations located at $w_2$ and $w_6$ to cover all routes. Given that we still need two stations even with reduced $D$ shows that our methodology is efficient in utilizing intersection segments to cover all routes with the minimum number of stations.

Table 7: Lists of endpoints when deviation distance $D = 3$ miles and $R = 100$ miles

| Endpoints $(w)$ | Set of covered routes $T(w)$ | Total flow covered $F(w)$ in round trips/day |
|---|---|---|
| $w_4, w_5, w_7$ | $U_1$ | 25 |
| $w_{36}, w_{55}$ | $U_2$ | 30 |
| $w_{44}, w_{45}$ | $U_3$ | 15 |
| $w_{13}, w_{18}$ | $U_4$ | 10 |
| $w_{49}, w_{57}, w_{58}, w_{59}$ | $U_5$ | 5 |
| $w_1, w_{12}$ | $U_1, U_2$ | 55 |
| $w_2$ | $U_1, U_2, U_3$ | 70 |
| $w_6$ | $U_1, U_4, U_5$ | 40 |
| $w_3, w_8, w_9, w_{10}, w_{11}, w_{48}, w_{56}$ | $U_1, U_5$ | 30 |
| $w_{25}, w_{26}, w_{27}, w_{31}, w_{32}, w_{41}, w_{42}$ | $U_2, U_3$ | 45 |
| $w_{22}, w_{23}, w_{33}$ | $U_2, U_4$ | 40 |
| $w_{24}, w_{30}, w_{34}, w_{52}, w_{53}$ | $U_2, U_4, U_5$ | 45 |
| $w_{35}, w_{37}, w_{38}, w_{39}, w_{40}, w_{54}$ | $U_2, U_5$ | 30 |
| $w_{15}, w_{16}, w_{19}, w_{20}, w_{21}, w_{29}, w_{46}, w_{47}$ | $U_3, U_4$ | 25 |
| $w_{43}$ | $U_3, U_5$ | 20 |
| $w_{14}, w_{17}, w_{28}, w_{50}, w_{51}$ | $U_4, U_5$ | 15 |

For Case 2, as the deviation distance increases, $D = 11$ miles, vehicles can now access more edges in the network, leading to a higher number of endpoints, 101 endpoints, and more intersection segments. The increase in endpoints suggests that a larger search space is created which can potentially lead to more alternative optimal stations and greater flexibility for the decision maker. Figure A.2 (see appendix) shows the set of endpoints for Case 2 with $D = 11$ miles and $R = 100$ miles and Table 8 shows the routes that can be covered by each endpoint.

Table 8: Lists of endpoints when deviation distance $D = 11$ miles and $R = 100$ miles

| Endpoints ($w$) | Set of covered routes $T(w)$ | Total flow covered $F(w)$ in round trips/day |
|---|---|---|
| $w_6$, $w_{36}$, $w_{39}$ | $U_1$, $U_2$, $U_3$ | 70 |
| $w_4$, $w_5$, $w_7$, $w_8$, $w_{11}$, $w_{12}$, $w_{13}$, $w_{14}$, $w_{15}$, $w_{19}$, $w_{20}$, $w_{21}$, $w_{24}$, $w_{25}$, $w_{26}$, $w_{27}$, $w_{28}$, $w_{29}$, $w_{31}$, $w_{41}$, $w_{42}$, $w_{45}$, $w_{52}$, $w_{68}$, $w_{88}$, $w_{89}$, $w_{90}$, $w_{98}$ | $U_1$, $U_2$, $U_3$, $U_4$ | 80 |
| $w_{43}$ | $U_1$,$U_2$, $U_3$, $U_4$, $U_5$ | 85 |
| $w_1$, $w_2$, $w_3$, $w_{35}$ | $U_1$, $U_2$, $U_3$, $U_5$ | 75 |
| $w_{34}$ | $U_1$, $U_2$, $U_5$ | 60 |
| $w_{30}$, $w_{92}$ | $U_1$, $U_3$, $U_4$ | 50 |
| $w_9$, $w_{10}$, $w_{16}$, $w_{17}$, $w_{18}$, $w_{32}$, $w_{71}$, $w_{72}$, $w_{73}$ | $U_1$, $U_3$, $U_4$, $U_5$ | 55 |
| $w_{22}$, $w_{23}$ | $U_1$, $U_3$, $U_5$ | 45 |
| $w_{33}$, $w_{86}$, $w_{91}$ | $U_1$, $U_4$, $U_5$ | 40 |
| $w_{38}$, $w_{49}$ | $U_2$, $U_3$ | 45 |
| $w_{44}$, $w_{46}$, $w_{47}$, $w_{48}$, $w_{50}$, $w_{51}$, $w_{53}$, $w_{54}$, $w_{55}$, $w_{56}$, $w_{63}$, $w_{67}$, $w_{80}$, $w_{87}$, $w_{93}$, $w_{94}$ | $U_2$, $U_3$, $U_4$ | 55 |
| $w_{40}$, $w_{58}$, $w_{62}$, $w_{64}$, $w_{65}$, $w_{74}$, $w_{75}$, $w_{76}$, $w_{77}$, $w_{78}$, $w_{81}$, $w_{82}$, $w_{83}$, $w_{95}$, $w_{95}$, $w_{99}$ | $U_2$, $U_3$, $U_4$, $U_5$ | 60 |
| $w_{66}$, $w_{100}$, $w_{101}$ | $U_2$, $U_3$, $U_5$ | 50 |
| $w_{57}$, $w_{59}$, $w_{60}$, $w_{61}$, $w_{79}$ | $U_2$, $U_4$, $U_5$ | 45 |
| $w_{37}$ | $U_2$, $U_5$ | 35 |
| $w_{69}$, $w_{70}$, $w_{84}$, $w_{85}$ | $U_3$, $U_4$, $U_5$ | 30 |
| $w_{96}$, $w_{97}$ | $U_4$, $U_5$ | 15 |

With an increased maximum allowed deviation distance, vehicles have greater flexibility to deviate from their dedicated routes to refuel. This allows vehicles operating on routes $U_1$, $U_2$, $U_3$, and $U_4$ to explore a wider range of refueling options. For example, a vehicle dedicated to route $U_2$ that deviates from vertex $v_{11}$ can now reach $w_{40}$, which lies on edge $(v_7, v_8)$. The increased flexibility in Case 2 can be seen in Table 8 where most endpoints can cover three or more routes. Notably, endpoint $w_{43}$ can cover the five routes, making it the optimal location for a single refueling station to cover all routes in the network. Note that route $U_5$ is an exception. The total length of route $U_5$ is 98 miles, and the vehicle's driving range $R$ is 100 miles. Consequently, to ensure that a vehicle dedicated to route $U_5$ can refuel, complete a full trip around its route, and come back to the refueling station, the maximum possible deviation distance will be limited by constraint (6) to 1 mile, even with a maximum allowed deviation distance $D = 11$ miles.

In Case 3, the driving range $R$ is increased from 100 to 106 miles while the deviation distance $D = 4$ miles. This increase in driving range affects the set of candidate refueling station locations by reducing the total number of endpoints to 56, as shown in Figure A.3 (see appendix) and Table 9.

In Case 2, vehicles dedicated to routes $U_1$, $U_2$, $U_3$, and $U_4$ could reach further edges compared to the original case because they are relatively short routes and the binding constraint was Constraint (7). Therefore, increasing $D$ helped them reach new edges in the network. However, route $U_5$ is long relative to the driving range $R$, therefore in the original case, the binding constraint was Constraint (6) which limited the possible deviation distance to 1 mile. In Case 3, the driving range $R$ is increased to 106 miles, now vehicles dedicated to route $U_5$ can deviate 4 miles instead of 1 mile before needing to refuel. In this case, two refueling stations are required, one can be located at $w_{35}$, $w_{49}$ or any point between them to cover routes $U_2$, $U_3$, $U_4$ and $U_5$, and a second station at $w_3$ to cover routes $U_1$ and $U_3$.

Table 9: Lists of endpoints when deviation distance $D = 4$ miles and $R = 106$ miles

| Endpoints ($w$) | Set of covered routes $T(w)$ | Total flow covered $F(w)$ in round trips/day |
|---|---|---|
| $w_{44}$ | $U_3$ | 15 |
| $w_2$, $w_{31}$ | $U_1$, $U_2$, $U_3$ | 70 |
| $w_1$, $w_{32}$ | $U_1$, $U_2$, $U_5$ | 60 |
| $w_3$ | $U_1$, $U_3$ | 40 |
| $w_4$, $w_{11}$ | $U_1$, $U_4$, $U_5$ | 40 |
| $w_5$, $w_6$, $w_7$, $w_8$, $w_9$, $w_{10}$ | $U_1$, $U_5$ | 30 |
| $w_{30}$, $w_{33}$, $w_{38}$, $w_{39}$, $w_{40}$, $w_{48}$, $w_{50}$, $w_{51}$, $w_{52}$, $w_{53}$ | $U_2$, $U_3$ | 45 |
| $w_{26}$ | $U_2$, $U_3$, $U_4$ | 55 |
| $w_{35}$, $w_{49}$ | $U_2$, $U_3$, $U_4$, $U_5$ | 60 |
| $w_{41}$, $w_{42}$ | $U_2$, $U_3$, $U_5$ | 50 |
| $w_{24}$, $w_{25}$, $w_{27}$, $w_{36}$ | $U_2$, $U_4$, $U_5$ | 45 |
| $w_{37}$, $w_{43}$ | $U_2$, $U_5$ | 35 |
| $w_{16}$, $w_{21}$, $w_{22}$, $w_{23}$, $w_{29}$, $w_{34}$, $w_{46}$ | $U_3$, $U_4$ | 25 |
| $w_{54}$, $w_{55}$, $w_{56}$ | $U_5$ | 5 |
| $w_{12}$, $w_{17}$, $w_{18}$, $w_{20}$, $w_{45}$, $w_{47}$ | $U_3$, $U_4$, $U_5$ | 30 |
| $w_{13}$, $w_{14}$, $w_{15}$, $w_{19}$, $w_{28}$ | $U_4$, $U_5$ | 15 |

This sensitivity analysis shows the effects of the two coverage constraints (6) and (7) and the values of parameters $D$ and $R$. In general, the higher the values of $D$ and $R$ the more flexibility vehicles have to deviate from their dedicated routes for refueling services. Increasing $D$ and $R$ allows vehicles to reach further locations in the network which increases the number and/or lengths of refueling segments (not necessarily the endpoints) providing more alternative optimal solutions. Moreover, since our methodology uses the two coverage constraints to find all refueling segments in the network and then extract a FDS, it remains robust under different parameter settings, consistently identifying the minimum number of stations while ensuring full coverage.

### *5.4. Computational Effort*

The ES algorithm is a polynomial-time algorithm with computational complexity $O\big(he(n^2 + h)\big)$ as shown in Theorem 3. Polynomial-time algorithms are generally regarded as computationally practical compared to exponential-time algorithms, which quickly become impractical for large-scale problems. In this section, we analyze the computational complexity of the proposed methodology for the two networks in Sections 5.1 and 5.2, thereby evaluating the scalability of the algorithm. After that, we compare the complexity of our methodology with that of well-known algorithms in the literature.

The Sioux Falls network has parameter values $h = 5, e = 40$, and $n = 27$. Based on Theorem 3, the number of operations in the dominant steps of the ES algorithm is 146,800. For the San Antonio network, the parameters increase to $h = 6, e = 52$, and $n = 32$, and the number of operations rose up to 321,360. Table 10 shows the computational times to implement the ES algorithm and solve the mathematical models using a PC with Windows 10 as the operating system and Intel core i9 3.00 GHz CPU with 8GB of Ram. The Sioux Falls network problem variations took between 0.014 and 0.036 seconds while the San Antonio network problem took 0.022 seconds.

Table 10: Computation time for different network scenarios

| Network | Case | No. of routes | Deviation distance ($D$) | Driving range ($R$) | No. of Endpoints | Computation time (seconds) |
|---|---|---|---|---|---|---|
| Sioux Falls | Original case | 5 | 4 | 100 | 67 | 0.022 |
| Sioux Falls | Case 1 | 5 | 3 | 100 | 59 | 0.014 |
| Sioux Falls | Case 2 | 5 | 11 | 100 | 101 | 0.036 |
| Sioux Falls | Case 3 | 5 | 4 | 106 | 56 | 0.014 |
| San Antonio | Original case | 6 | 4 | 100 | 77 | 0.022 |

Figure 15 shows the growth in number of operations in the dominant steps of the ES algorithm for various numbers of vertices and routes, while assuming the number of edges in a complete graph, $e = n(n-1)/2$. We plotted the number of operations for $h$ values between 5 and 350 routes. To put things in perspective, San Antonio, TX is the 7$^{\text{th}}$ largest city in the USA by population [52]. VIA Metropolitan Transit serves San Antonio with only 75 bus routes [53]. Moreover, the Metropolitan Transportation Authority is North America's largest transportation network serving New York City [54], has about 327 bus routes. As can be seen in the figure, even with a large number of routes, $h = 350$, the computational complexity of the ES algorithm remains close to practical polynomial complexities such as $O(n^4)$ and $O(n^5)$. These complexities are observed in well-known and widely used methods. For instance, Edmonds and Karp [55] algorithm for solving the maximum flow problem has the complexity of $O(n^5)$. Another example is the longest path problem, as described in [56], has the complexity of $O(n^4)$. Conversely, exponential complexities, are often encountered in brute-force solutions to combinatorial optimization problems. For instance, the brute-force algorithm for solving the Knapsack problem operates in $O(2^n)$ time, meaning the complexity grows exponentially as the problem size increases, making it computationally intractable [57]. Therefore, our proposed algorithm is practical even for realistic problem instances.

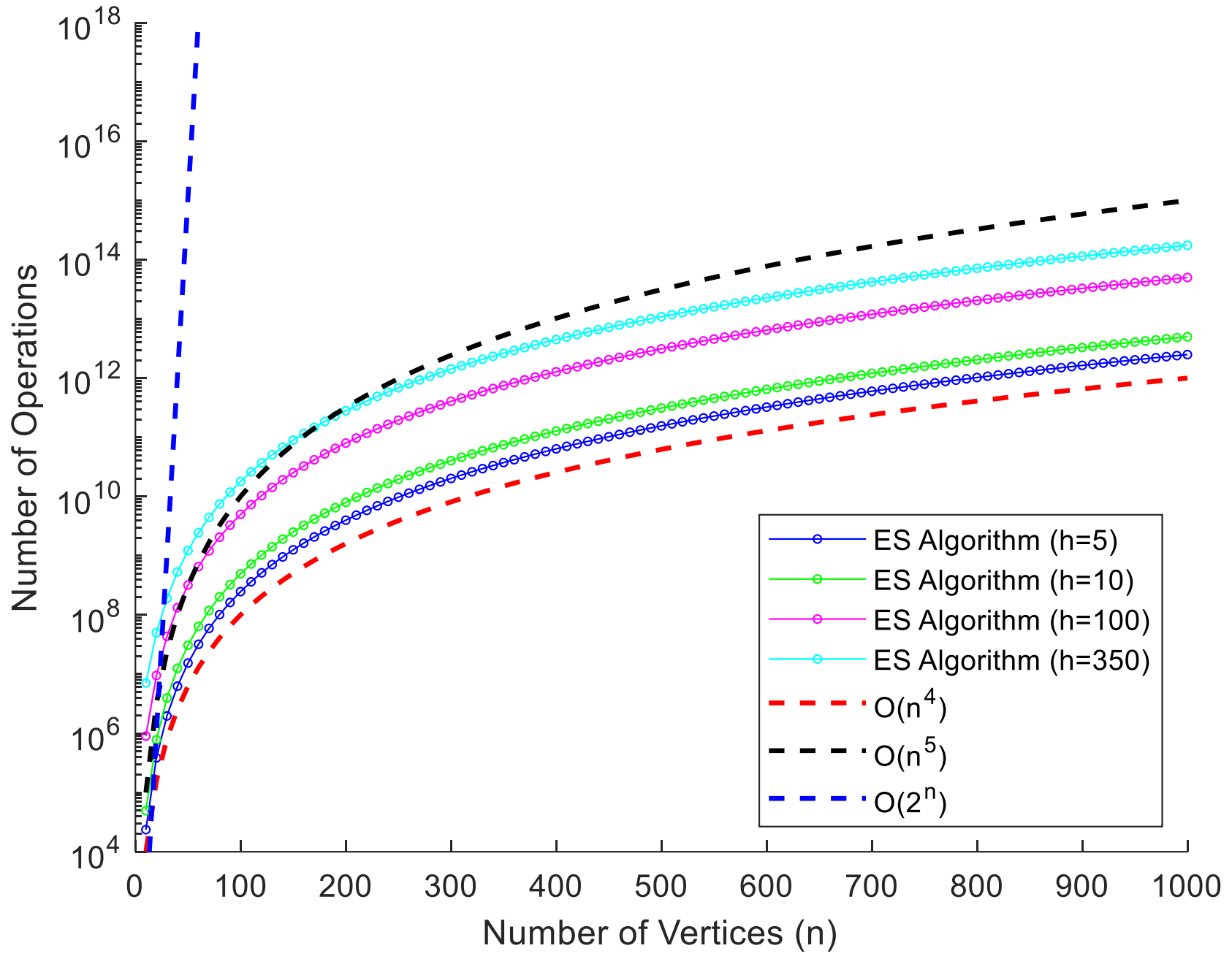

Figure 15: Computational complexity growth of the ES algorithm

## 6. Conclusions and Future Work

In this paper, we address the refueling station location problem for dedicated fleet vehicles with pre-planned fixed routes. Most existing literature in this area follows a discrete approach that assumes a predetermined set of candidate locations and then finds a solution among them. The discrete approach cannot guarantee a global optimal solution. To the best of our knowledge, this is the first paper following the continuous approach to solve the refueling station location problem for dedicated routes on a general network. In this approach, we assume that all points in the network are candidate locations, i.e., infinite set of candidate locations. We develop the mathematical foundation for an exact algorithm to extract a FDS guaranteed to have a global optimal solution. Therefore, combining optimality and computational efficiency.

The rise of AF-powered vehicles characterized by limited driving ranges made it necessary for route planners and service providers to ensure that these vehicles can serve their customers without running out of fuel. On one hand, the underdeveloped and sparse AF refueling infrastructure force drivers to deviate from their pre-planned routes for refueling services. On the other hand, passenger convenience, strict arrival and departure time windows, and narrow operating ratios for carriers restrict the acceptable deviation distance a vehicle can drive for refueling services. In this study, we include these technical and operational constraints to develop a practical approach that can be applied in realistic problems.

The proposed methodology scans all network edges for potential refueling sets that can cover the dedicated routes in the network. For a refueling set to cover a route it must satisfy two constraints. First, vehicles of the covered route must be able to travel at least one trip around their dedicated route and then reach any point in the refueling set before running out of fuel. Second, vehicles must not exceed a predetermined maximum allowed deviation distance from their dedicated route to reach a refueling station. We prove that the endpoints of these refueling sets constitute a FDS for our problem. The algorithm can find the overall set of endpoints in $O\big(he(n^2 + h)\big)$. After that, the set of endpoints can be used in a set covering model to find the minimum number of refueling stations required to cover all routes in the network and their locations. Alternative optimal solutions can be found by replacing endpoints selected in the optimal solution by other endpoints or interior points covering similar routes. We prove that all interior points of a line segment connecting two adjacent endpoints cover the same set of routes. Therefore, if an alternative optimal solution can be obtained by replacing an endpoint by an interior point of a line segment between two adjacent endpoints, we can obtain an infinite number of alternative optimal solutions. This flexibility allows decision-makers to select an alternative optimal solution that fits their needs considering other criteria like land prices, accessibility, slope, etc.

Sensitivity analysis reveals that the maximum allowed deviation distance and driving range play crucial roles in determining the number of refueling stations required to cover the flows in a network. Smaller deviation allowances and shorter ranges restrict the number of accessible edges outside a dedicated route, leading to fewer and shorter intersection refueling segments where a single refueling station can serve multiple routes. In contrast, allowing longer deviation distances and greater driving ranges enables vehicles to reach more edges, thereby reducing the number of required refueling stations. Both theoretical results and numerical experiments demonstrate that the proposed methodology is efficient and practical, exhibiting computational complexity comparable to many well-known and widely used algorithms.

Potential future work ideas include considering multiple vehicle classes with different driving ranges. Also, our assumption of a predetermined and fixed maximum allowed deviation distance can be replaced by a percentage of the route length. This change will allow vehicles with long routes to deviate more toward a refueling station while restricting vehicles with short routes to deviate less. Finally, the current business practice is that each company has its own AF refueling station, usually at a larger service shop. However, this option may be infeasible for smaller fleets. A potential research study can focus on designing a collaborative structure between smaller fleets that do not have dedicated refueling stations and larger fleets.

## APPENDIX A

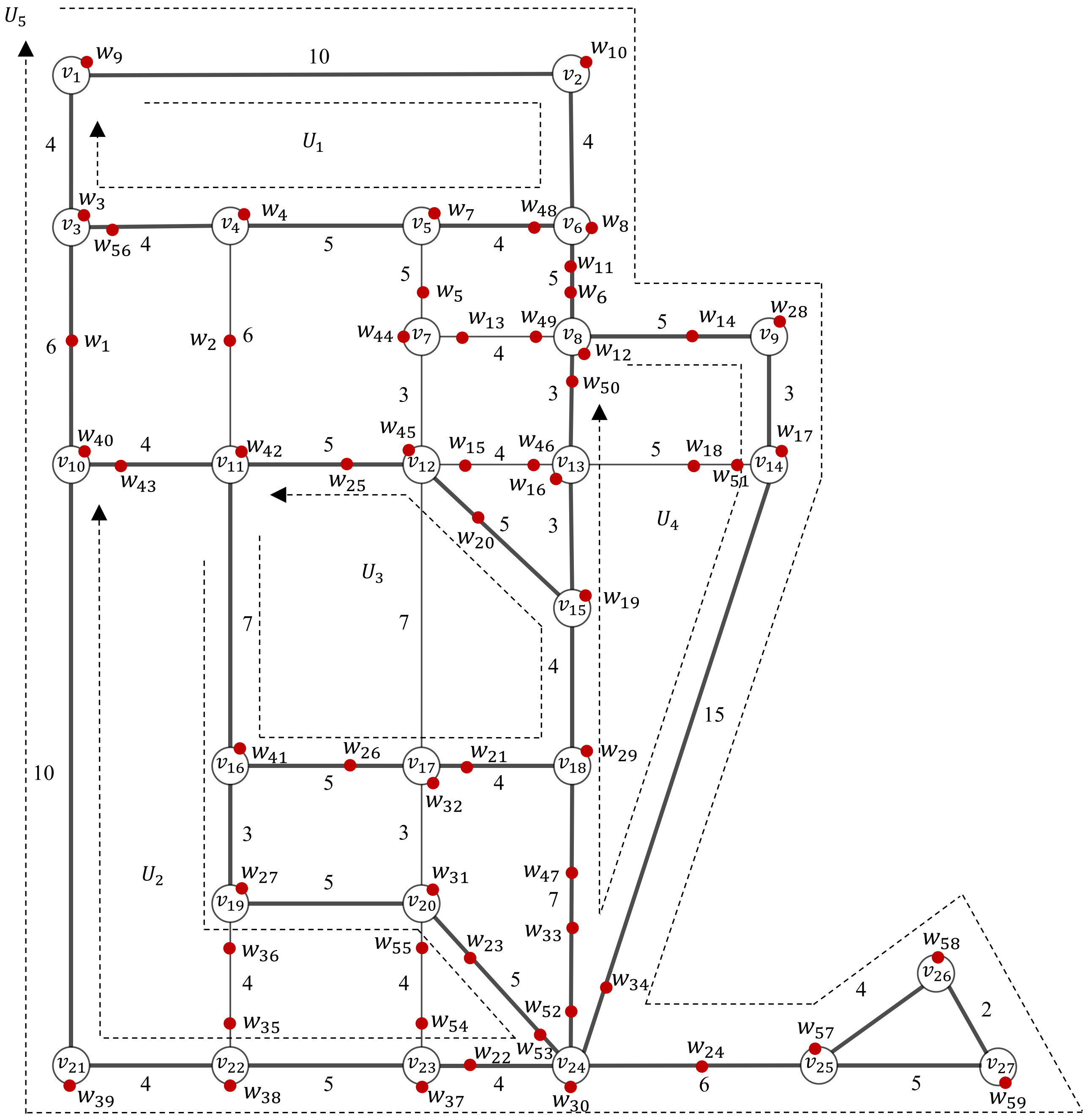

Figure A.1: Dedicated route example when D=3 miles and R=100 miles

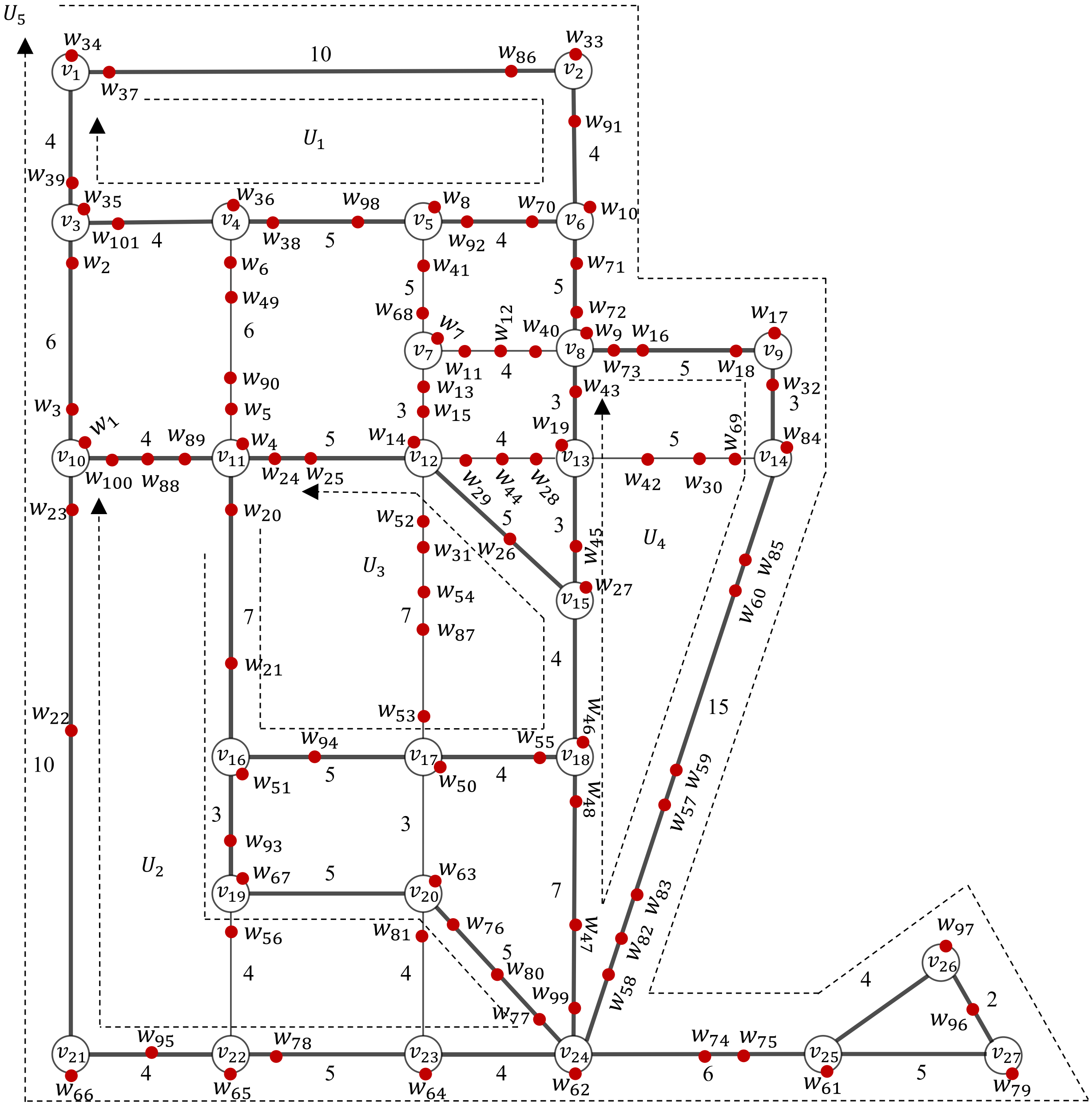

Figure A.2: Dedicated route example when D=11 miles and R=100 miles

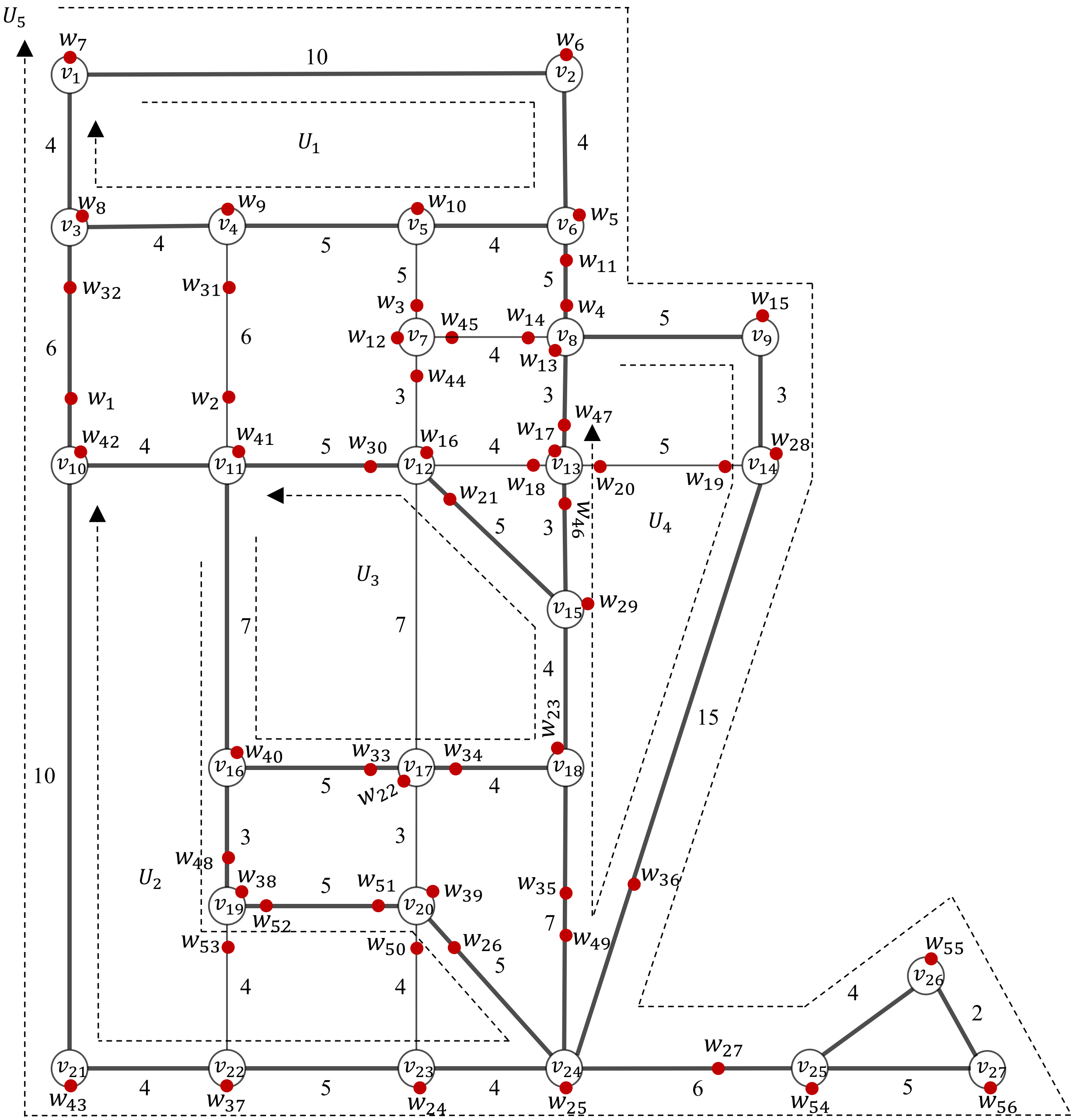

Figure A.3: Dedicated route example when D=4 miles and R=106 miles